\documentclass[a4paper,12pt,twoside]{article}
\usepackage[british]{babel}
\usepackage{amssymb,amsmath,amsfonts}
\usepackage{graphicx}
\usepackage{subfigure}
\newtheorem{theorem}{Theorem}[section]

\def\vet#1{{\boldsymbol{#1}}}
\def\build#1_#2^#3{\mathrel{
\mathop{\kern 0pt#1}\limits_{#2}^{#3}}}
\def\reali{\mathbb{R}}
\def\complessi{\mathbb{C}}

\def\interi{\mathbb{Z}}
\def\toro{\mathbb{T}}
\def\Ascr{\mathcal{A}}
\def\Bscr{\mathcal{B}}
\def\Cscr{\mathcal{C}}
\def\Dscr{\mathcal{D}}
\def\Escr{\mathcal{E}}

\def\Gscr{\mathcal{G}}
\def\Hscr{\mathcal{H}}
\def\Kscr{\mathcal{K}}
\def\Lscr{\mathcal{L}}
\def\Oscr{\mathcal{O}}
\def\Pscr{\mathcal{P}}
\def\Sscr{\mathcal{S}}
\def\Tscr{\mathcal{T}}
\def\lie#1{\Lscr_{#1}}

\def\epsilon{\varepsilon}
\def\rho{\varrho}
\def\phi{\varphi}
\def\csi{\xi}
\def\imunit{{\rm i}}
\def\poisson#1#2{\lbrace #1,#2 \rbrace}
\def\adaptpoisson#1#2{\left\{ #1\,,\,#2 \right\}}

\def\scalprod#1#2{#1\cdot#2}

\title{\bf A Semi-Analytic Algorithm for Constructing Lower Dimensional
  Elliptic Tori in Planetary Systems.
  \thanks{ 
    Mechanics.  {\it 2010 Mathematics Subject Classification.}
    Primary: 37J40; Secondary: 37N05, 70F10, 70--08, 70H08.  }}

\author{
  {\bf MARCO SANSOTTERA}\\
  {\small Dipartimento di Matematica, Universit\`a degli Studi di Milano,}\\
  {\small via Saldini 50, 20133\ ---\ Milano, Italy.}\\
  {\bf UGO LOCATELLI}\\
  {\small Dipartimento di Matematica, 
    Universit\`a degli Studi di Roma ``Tor Vergata'',}\\
  {\small Via della Ricerca Scientifica 1, 00133--Roma (Italy).}\\
  {\bf ANTONIO GIORGILLI}\\
  {\small Dipartimento di Matematica, Universit\`a degli Studi di Milano,}\\
  {\small via Saldini 50, 20133\ ---\ Milano, Italy.}\\
  {\small e-mails:
    {\tt marco.sansottera@gmail.com,$\>$locatell@mat.uniroma2.it}}\\
  {\small {\tt antonio.giorgilli@unimi.it}}
}

\date{}

\begin{document}
\maketitle

\begin{abstract}
We adapt the Kolmogorov's normalization algorithm (which is the key
element of the original proof scheme of the KAM theorem) to the
construction of a suitable normal form related to an invariant
elliptic torus.  As a byproduct, our procedure can also provide some
analytic expansions of the motions on elliptic tori. By extensively
using algebraic manipulations on a computer, we explicitly apply our
method to a planar four-body model not too different with respect to
the real Sun--Jupiter--Saturn--Uranus system. The frequency analysis
method allows us to check that our location of the initial conditions
on an invariant elliptic torus is really accurate. 
\end{abstract}

\bigskip

\section{Introduction}\label{sec:intro}

Since the birth of KAM theory (see~\cite{Kolmogorov-1954},
\cite{Moser-1962} and~\cite{Arnold-1963.1}), invariant tori are
expected to be the key dynamical objects which may explain the (nearly
perfect) quasi-periodicity of the planetary motions of our Solar
System.  On the other hand, Poincar\'e has widely discussed the role
of periodic orbits, namely of one-dimensional tori.  A natural
extension is the search for lower dimensional invariant tori, which
are characterized by quasi-periodic motions with a number of
frequencies lower than the actual number of degrees of freedom of the
system.

In the case of a planetary system including $n$ planets it seems
natural to look for particular orbits, namely solutions of Newton's
equations lying on $n$-dimensional tori, the latter being slight
deformations of the composition of $n$ {\it coplanar circular\/}
Keplerian orbits.  This in view of Lagrange theory of secular motions,
where, in the frame of the so-called {\it secular model\/}, circular
orbits are considered as a first approximation and the motions of
eccentricities and inclinations are represented as small oscillations
with respect to the reference orbit.  Replacing the circular orbit
with an elliptic lower dimensional torus, i.e., a
torus which is stable with respect to small oscillations around it,
we expect to find a better starting point for a secular theory.

The aim of this paper is to present an explicit algorithm, based on
the original Kolmogorov's one, for constructing elliptic tori for a
model including three of the biggest planets.  We also perform an
explicit calculation via algebraic manipulations and a comparison of
the orbits so found with the outcome of a numerical integration of
Newton's equations.

It should be emphasized that the existence of elliptic lower
dimensional invariant tori is not a straightforward consequence of
Kolmogorov's theorem.

The existence of lower dimensional tori in the planetary problem (both
elliptic or hyperbolic) has been proven by Jefferys and Moser
(see~\cite{Jef-Mos-1966}) and Lieberman (see~\cite{Lieberman-1971}).
A general theorem has been proved by P{\"o}schel
(see~\cite{Poschel-1989} and~\cite{Poschel-1996}).  An application of
P{\"o}schel's method to the Solar System has been produced by Biasco,
Chierchia and Valdinoci in two different cases, namely the spatial
three-body planetary problem and a planar system with a central
star and $n$ planets (see~\cite{Bia-Chi-Val-2003}
and~\cite{Bia-Chi-Val-2006}, respectively).  However, as it often
happens in the framework of KAM theory, their approach is a deep one
from a theoretical point of view, but seems not to be suitable for
explicit calculations, even if one is interested just in finding the
locations of the elliptic invariant tori.  Let us clarify this point.

The results of Biasco, Chierchia and Valdinoci are based on performing
a sequence of canonical transformations that give the Hamiltonian a
particular form to which P{\"o}schel's theorem can be applied.  In
turn, P{\"o}schel's theorem is based on a clever adaptation of
Arnold's proof of KAM theorem.  Precisely, the perturbation is removed
by a sequence of canonical transformations which are defined on a
subset of the phase space excluding the ``resonant regions''
(see~\cite{Arnold-1963.1} and~\cite{Arnold-1963.2}).  The procedure
ends up with a nowhere dense Cantor set of elliptic invariant tori.
Although effective from the analytical viewpoint, the procedure looks
hardly applicable in a practical calculation.

The original scheme of the proof introduced by Kolmogorov is in a much
better position for what concerns the translation into an explicit
algorithm (see~\cite{Kolmogorov-1954}, \cite{Ben-Gal-Gio-Str-1984},
\cite{Gio-Loc-1997.1} and~\cite{Gio-Loc-1997.2}).  Such an approach
has been successfully used to calculate the orbits for some
interesting problems in Celestial Mechanics (see~\cite{Loc-Gio-2000},
\cite{Loc-Gio-2005}, \cite{Loc-Gio-2007} and~\cite{Gab-Jor-Loc-2005}).
Thus, we think it is useful to modify Kolmogorov's algorithm so as to
be able to construct a suitable normal form related to the elliptic
tori.  In addition, this will allow us to explicitly integrate the
equations of motion on those invariant surfaces, by using a so-called
semi-analytic procedure.

An explicit construction of elliptic (or maybe hyperbolic) lower
dimensional tori may be useful in many cases.  We pay some attention
to (a)~the problem of the long term stability of orbits in the
planetary system, and (b)~the possible application to the search for
orbits suitable for spatial missions.

Concerning the long term stability of a planetary system, the
construction of a normal form related to some fixed elliptic torus
could be a relevant milestone.  It is indeed possible to ensure the
effective stability in the neighborhood of such an invariant surface
by implementing a partial construction of the Birkhoff normal form.  A
similar approach in the neighborhood of either a maximal dimensional
KAM torus or of the circular orbits has been worked out, e.g.,
in~\cite{Jor-Vil-1997.2}, \cite{Gio-Loc-San-2009} and~\cite{Gio-Loc-San-2010}.  Concerning our
Solar System, such an approach might be applied to some asteroids with
small orbital eccentricities and inclinations.  However, as explained
in~\cite{San-Loc-Gio-2010}, this same approach can not yet succeed in
proving the long-time stability of the major planets of our Solar
System.

Concerning spatial missions, the strong stability of the regions close
to elliptic invariant tori may be very useful in identifying stable
orbits.  Moreover, our technique should also adapt quite easily to the
construction of hyperbolic tori that can be used in the design of
spacecraft missions requiring low energy transfers.  We also recall
that lower dimensional tori of elliptic, hyperbolic and mixed type
have been studied in the vicinity of the Lagrangian points for both
the restricted three-body problem and the bicircular restricted
four-body problem (see, e.g.,~\cite{Jor-Vil-1997.1},
\cite{Jor-Vil-1998}, \cite{Cas-Jor-2000} and~\cite{Gab-Jor-2001}).

The paper is devoted to the construction of elliptic tori for a model
not far from the Sun--Jupiter--Saturn--Uranus (SJSU) {\it planar\/}
system (let us recall that the real orbits of the planets of our Solar
System {\it are presumably not\/} lying on lower dimensional tori).  By the way, we
think that with some minor modifications our procedure should adapt
also to the more general spatial case, after having performed the
reduction of the angular momentum, which is not considered here in
order to shorten the description of all the preliminary expansions
(for an introduction to some methods performing both the partial and
the total reduction, see~\cite{Deprit-1983}, \cite{Mal-Rob-Las-02}
and~\cite{Pinzari-tesi}).  The contents are organized as follows.

Sect.~\ref{sec:2D_plan_Ham} is devoted to the introduction of our
Hamiltonian model and to a schematic description of its expansion in
canonical coordinates. This will allow us to write down the form of
the Hamiltonian to which our approach can be applied.

Our algorithm constructing a normal form for elliptic tori is
presented in a purely formal way in sect.~\ref{sec:form_alg}.  In
subsect.~\ref{sbs:remarks_theory} we briefly introduce the theoretical
background necessary to investigate the convergence.  We plan to give
all details about convergence in a future work
(see~\cite{Gio-Loc-San-2011}).

Sect.~\ref{sec:num_tests} is devoted to an application which we also
use as a test of our procedure.  First, in
subsects.~\ref{sbs:comp_algebra}--\ref{sbs:expl_calc} we describe the
implementation of our algorithm, by using an algebraic manipulator on
a computer so to produce both the normal form and the semi-analytic
integration of the motion on an invariant elliptic torus.  Then in
sect.~\ref{sbs:val_freq_anal} we check the accuracy of our
construction by using frequency analysis.  This is based on the fact
that the Fourier spectrum of the motions on elliptic tori is strongly
characteristic, because only the mean-motion frequencies and their
linear combinations can show up.

\section{Classical Expansion of the Planar
  Planetary Hamiltonian}\label{sec:2D_plan_Ham}

We apply our algorithm for constructing elliptic tori to a concrete
planetary model.  We consider four point bodies
$P_0,\,P_1,\,P_2,\,P_3$, with masses $m_0,\,m_1,\,m_2,\,m_3$, mutually
interacting according to Newton's gravitational law.  Hereafter we
associate the indexes $0,\,1,\,2,\,3$ to Sun, Jupiter, Saturn and
Uranus, respectively.\footnote{Let us stress that the four bodies have
  the same masses as Sun, Jupiter, Saturn and Uranus, but the orbits
  studied here are significantly different with respect to the real
  ones.}

Let us now recall how the classical Poincar\'e's variables can be
introduced so to perform a first expansion of the Hamiltonian around
circular orbits, i.e., having zero eccentricity.  We basically follow
the formalism introduced by Poincar\'e (see \cite{Poincare-1892}
and~\cite{Poincare-1905}; for a modern exposition, see,
e.g.,~\cite{Laskar-1989b} and~\cite{Las-Rob-1995}).  We remove the
motion of the center of mass by using heliocentric coordinates
$\vet{r}_j=\build{P_0P_j}_{}^{\longrightarrow}\,$, with
$j=1,\,2,\,3\,$.  Denoting by $\tilde{\vet{r}}_j$ the momenta
conjugated to $\vet{r}_j$, the Hamiltonian of the system has $6$
degrees of freedom, and reads
\begin{equation}
F(\tilde{\vet{r}},\vet{r})= T^{(0)}(\tilde{\vet{r}})+U^{(0)}(\vet{r})+
T^{(1)}(\tilde{\vet{r}})+U^{(1)}(\vet{r}) \,,
\label{Ham-iniz}
\end{equation}
where
$$
\begin{array}{rclrcl}
T^{(0)}(\tilde{\vet{r}}) &= &\frac{1}{2}\build{\sum}_{j=1}^{3} 
\frac{m_0+m_j}{m_0m_j}\,\|\tilde{\vet{r}}_j\|^2\,,
\qquad
&T^{(1)}(\tilde{\vet{r}}) &=
&\frac{1}{m_0}\Big(\tilde{\vet{r}}_1\cdot\tilde{\vet{r}}_2+
\tilde{\vet{r}}_1\cdot\tilde{\vet{r}}_3+
\tilde{\vet{r}}_2\cdot\tilde{\vet{r}}_3\Big)\,,
\cr\cr
U^{(0)}(\vet{r}) &= &-\Gscr\build{\sum}_{j=1}^{3}
\frac{m_0\, m_j}{\|\vet{r}_j\|}\,,
\qquad
&U^{(1)}(\vet{r})
&= &-\Gscr\left(\frac{m_1\, m_2}{\|\vet{r}_1-\vet{r}_2\|}
+\frac{m_1\, m_3}{\|\vet{r}_1-\vet{r}_3\|}
+\frac{m_2\, m_3}{\|\vet{r}_2-\vet{r}_3\|}\right)\,.
\cr
\end{array}
$$

The plane set of Poincar\'e's canonical variables is introduced as
\begin{equation}
\vcenter{\openup1\jot\halign{
 \hbox {\hfil $\displaystyle {#}$}
&\hbox {\hfil $\displaystyle {#}$\hfil}
&\hbox {$\displaystyle {#}$\hfil}
&\hbox to 6 ex{\hfil$\displaystyle {#}$\hfil}
&\hbox {\hfil $\displaystyle {#}$}
&\hbox {\hfil $\displaystyle {#}$\hfil}
&\hbox {$\displaystyle {#}$\hfil}\cr
\Lambda_j &=& \frac{m_0\, m_j}{m_0+m_j}\sqrt{\Gscr(m_0+m_j) a_j}\,,
& &\lambda_j &=& M_j+\omega_j\,,
\cr
\csi_j &=& \sqrt{2\Lambda_j}
\sqrt{1-\sqrt{1-e_j^2}}\,\cos\omega_j\,,
& &\eta_j&=&-\sqrt{2\Lambda_j}
\sqrt{1-\sqrt{1-e_j^2}}\, \sin\omega_j\,,
\cr
}}
\label{var-Poincare-piano}
\end{equation}
for $j=1\,,\,2\,,\,3\,$, where $a_j\,,\> e_j\,,\> M_j$ and $\omega_j$
are the semi-major axis, the eccentricity, the mean anomaly and the
perihelion argument, respectively, of the $j$-th planet.  One
immediately sees that both $\csi_j$ and $\eta_j$ are of the same order
of magnitude as the eccentricity $e_j\,$. Using Poincar\'e's
variables~(\ref{var-Poincare-piano}), the Hamiltonian $F$ can be
rearranged so that one has
\begin{equation}
F(\vet{\Lambda},\vet{\lambda},\vet{\csi},\vet{\eta})=
F^{(0)}(\vet{\Lambda})+
F^{(1)}(\vet{\Lambda},\vet{\lambda},\vet{\csi},\vet{\eta}) \,,
\label{Ham-iniz-Poincare-var}
\end{equation}
where $F^{(0)}=T^{(0)}+U^{(0)}$, $F^{(1)}=T^{(1)}+U^{(1)}$.  Let us
emphasize that $F^{(0)}=\Oscr(1)$ and $F^{(1)}=\Oscr(\mu)\,$, where
the small dimensionless parameter
$\mu=\max\{m_1\,/\,m_0\,,\,m_2\,/\,m_0\,,\,m_3\,/\,m_0\,\}$ highlights
the different size of the terms appearing in the Hamiltonian.
Therefore, let us remark that the time derivative of each coordinate
is $\Oscr(\mu)$ but in the case of the angles $\vet{\lambda}\,$. Thus,
according to the common language in Celestial Mechanics, in the
following we will refer to $\vet{\lambda}$ and to their conjugate
actions $\vet{\Lambda}$ as the {\em fast variables}, while
$(\vet{\csi},\vet{\eta})$ will be called {\em secular variables}.

We proceed now by expanding the
Hamiltonian~(\ref{Ham-iniz-Poincare-var}) in order to construct the
first basic approximation of the normal form for elliptic tori.  After
having chosen a center value $\vet{\Lambda}^*$ for the Taylor
expansions with respect to the fast actions (in a way we will explain
later), we perform a translation $\Tscr_{\vet{\Lambda}^*}$ defined as
\begin{equation}
L_j=\Lambda_j-\Lambda_j^*\,,
\qquad\forall\> j=1\,,\, 2\,,\, 3\,.
\label{def-L}
\end{equation}
This is a canonical transformation that leaves the coordinates
$\vet{\lambda}\,$, $\vet{\csi}$ and $\vet{\eta}$ unchanged.  The
transformed Hamiltonian
$\Hscr^{(\Tscr)}=F\circ\Tscr_{\vet{\Lambda}^*}\,$ can be expanded in
power series of $\vet{L},\,\vet{\csi},\,\vet{\eta}$ around the origin.
Thus, forgetting an unessential constant we rearrange the Hamiltonian
of the system as
\begin{equation}
\Hscr^{(\Tscr)}(\vet{L},\vet{\lambda},\vet{\csi},\vet{\eta})=
\vet{n}^{*}\cdot\vet{L}+
\sum_{j_1=2}^{\infty}h_{j_1\,,\,0}^{({\rm Kep})}(\vet{L})+
\sum_{j_1=0}^{\infty}\sum_{j_2=0}^{\infty}
h_{j_1\,,\,j_2}^{(\Tscr)}(\vet{L},\vet{\lambda},\vet{\csi},\vet{\eta}) \,,
\label{Ham-trasl-fast}
\end{equation}
where the functions $h_{j_1\,,\,j_2}^{(\Tscr)}$ are homogeneous
polynomials of degree $j_1$ in the actions $\vet{L}$ and of degree
$j_2$ in the secular variables $(\vet{\csi},\vet{\eta})\,$. The
coefficients of such homogeneous polynomials do depend analytically
and periodically on the angles $\vet{\lambda}\,$.  The terms
$h_{j_1\,,\,0}^{({\rm Kep})}$ of the Keplerian part are homogeneous
polynomials of degree $j_1$ in the actions $\vet{L}\,$, the explicit
expression of which can be determined in a straightforward manner.  In
the latter equation the term which is both linear in the actions and
independent of all the other canonical variables (i.e.,
$\vet{n}^{*}\cdot\vet{L}$) has been separated in view of its relevance
in perturbation theory, as it will be discussed in the next section.
We also expand the coefficients of the power series
$h_{j_1\,,\,j_2}^{(\Tscr)}$ in Fourier series of the angles
$\vet{\lambda}\,$.  The expansion of the Hamiltonian is a traditional
procedure in Celestial Mechanics.  We work out these expansions for
the case of the planar SJSU system using a specially devised algebraic
manipulator.  The calculation is based on the approach described in
sect.~2.1 of~\cite{Loc-Gio-2000}, which in turn uses the scheme
sketched in sect.~3.3 of~\cite{Robutel-1995}.

\begin{table*}
\caption[]{Masses $m_j$ and initial conditions for Jupiter, Saturn and
  Uranus in our planar model.  We adopt the AU as unit of length, the
  year as time unit and set the gravitational constant
  $\Gscr=1\,$. With these units, the solar mass is equal to
  $(2\pi)^2$. The initial conditions are expressed by the usual
  heliocentric planar orbital elements: the semi-major axis $a_j\,$,
  the mean anomaly $M_j\,$, the eccentricity $e_j$ and the perihelion
  longitude $\omega_j\,$.  The data are taken by JPL at the Julian
  Date $2440400.5\,$.  }
\label{tab:parameters_2D_SJSU}
\begin{center}
  \begin{tabular}{|c|l|l|l|}
\hline
& Jupiter ($j=1$) & Saturn ($j=2$) & Uranus ($j=3$)
\\
\hline
$m_{j}^{\phantom{\displaystyle 1}}$
& $(2\pi)^2/1047.355$
& $(2\pi)^2/3498.5$
& $(2\pi)^2/22902.98$
\\
$a_j$
& $5.20463727204700266$ & $9.54108529142232165$ & $19.2231635458410572$
\\
$M_j$
& $3.04525729444853654$ & $5.32199311882584869$ & $0.19431922829271914$
\\
$e_j$
& $0.04785365972484999$ & $0.05460848595674678$ & $0.04858667407651962$
\\
$\omega_j$
& $0.24927354029554571$ & $1.61225062288036902$ & $2.99374344439246487$ 
\\
\hline
\end{tabular}
\end{center}
\end{table*}

The reduction to the planar case is performed as follows.  We pick
from table~IV of~\cite{Standish-1998} the initial conditions of the
planets in terms of heliocentric positions and velocities at the
Julian Date $2440400.5\,$.  Next, we calculate the corresponding
orbital elements with respect to the invariant plane (that is
perpendicular to the total angular momentum).  Finally we include the
longitudes of the nodes $\Omega_j$ (which are meaningless in the
planar case) in the corresponding perihelion longitude $\omega_j$ and
we eliminate the inclinations by setting them equal to zero.  The
remaining initial values of the orbital elements are reported in
table~\ref{tab:parameters_2D_SJSU}.

Having fixed the initial conditions we come to determining the {\it
  average\/} values $(a_1^*, a_2^*, a_3^*)$ of the semi-major axes
during the evolution.  To this end we perform a long-term numerical
integration of Newton's equations starting from the initial conditions
related to the data reported in table~\ref{tab:parameters_2D_SJSU}.
After having computed $(a_1^*, a_2^*, a_3^*)\,$, we determine the
values $\vet{\Lambda}^*$ via the first equation
in~(\ref{var-Poincare-piano}).  This allows us to perform the
expansion~(\ref{Ham-trasl-fast}) of the Hamiltonian as a function of
the canonical coordinates
$(\vet{L},\vet{\lambda},\vet{\csi},\vet{\eta})$.  In our calculations
we truncate this initial expansion as follows.  (a)~The Keplerian part
is expanded up to the quartic terms.  The series where the general
summand $h_{j_1\,,\,j_2}^{(\Tscr)}$ appears are truncated so to
include: (b1)~all terms having degree $j_1$ in the actions $\vet{L}$
with $j_1\le 3\,$, (b2)~all terms having degree~$j_2$ in the secular
variables $(\vet{\csi},\vet{\eta})\,$, with $j_2$ such that
$2j_1+j_2\le 8\,$, (b3)~all terms up to the trigonometric degree $18$
with respect to the angles $\vet{\lambda}\,$.  This choice is
motivated by the remark that the orbits on elliptic tori reach values
of the eccentricities smaller than those attained by the real motions
(let us recall that both $\csi_j=\Oscr(e_j)$ and $\eta_j=\Oscr(e_j)$
$\forall\ j=1\,,\,2\,,\,3\,$); moreover, larger limits on the fast
angles are needed, in order to give a sharp enough numerical
evidence of the convergence of the algorithm described in the next
section.

Let us now focus on the average with respect to the fast angles of the
Hamiltonian written in~(\ref{Ham-trasl-fast}), i.e.,
$\big\langle\Hscr^{(\Tscr)}\big\rangle_{\vet{\lambda}}\,$.  The fast
actions $\vet{L}$ are obviously invariant with respect to the flow of
$\big\langle\Hscr^{(\Tscr)}\big\rangle_{\vet{\lambda}}\,$, thus, they
can be neglected inasmuch as the secular motions are considered. The
remaining most significant term is given by the lowest order
approximation of the secular Hamiltonian, namely its quadratic part
$\big\langle h_{0,2}^{(\Tscr)}\big\rangle_{\vet{\lambda}}$ which is
essentially the one considered in the theory first developed by
Lagrange (see~\cite{Lagrange-1776}) and further improved by Laplace
(see~\cite{Laplace-1772}, \cite{Laplace-1784} and~\cite{Laplace-1785})
and by Lagrange himself (see~\cite{Lagrange-1781},
\cite{Lagrange-1782}).  In modern language, we can say that the origin
$(\vet{\csi},\vet{\eta})=(\vet{0},\vet{0})\,$ of the secular
coordinates is an elliptic equilibrium point for the secular
Hamiltonian.  It is well known that under mild assumptions on the
quadratic part of the Hamiltonian which are satisfied in our case (see
sect.~3 of~\cite{Bia-Chi-Val-2006} where such hypotheses are shown to
be generically fulfilled for a planar model of our Solar System) one
can find a canonical transformation
$(\vet{L},\vet{\lambda},\vet{\csi},\vet{\eta})=
\Dscr(\vet{p},\vet{q},\vet{x},\vet{y})$ with the following properties:
(i)~$\vet{L}=\vet{p}$ and $\vet{\lambda}=\vet{q}\,$, (ii)~the map
$(\vet{\csi},\vet{\eta})=
\big(\vet{\csi}(\vet{x}),\vet{\eta}(\vet{y})\big)$ is linear,
(iii)~$\Dscr$ diagonalizes the quadratic part of the Hamiltonian, so
that we can write $\big\langle
h_{0,2}^{(\Tscr)}\big\rangle_{\vet{\lambda}}$ in the new coordinates
as $\sum_{j=1}^3\nu_j^{(0)}(x_j^2+y_j^2)/2\,$, where all the entries
of the vector $\vet{\nu}^{(0)}$ have the same sign.  Our algorithm
constructing a suitable normal form for elliptic tori starts from the
Hamiltonian $H^{(0)}=\Hscr^{(\Tscr)}\circ\Dscr\,$, i.e.,
\begin{equation}
H^{(0)}(\vet{p},\vet{q},\vet{x},\vet{y})=
\Hscr^{(\Tscr)}\big(\Dscr(\vet{p},\vet{q},\vet{x},\vet{y})\big)\,.
\label{Ham^(0)}
\end{equation}

\section{Formal Algorithm}\label{sec:form_alg}

We present the formal algorithm making reference to a generic
Hamiltonian with $n_1+n_2\,$ degrees of freedom, where the canonical
coordinates $(\vet{p},\vet{q},\vet{x},\vet{y})$ can naturally be split
in two parts, that are
$(\vet{p},\vet{q})\in\reali^{n_1}\times\toro^{n_1}$ and
$(\vet{x},\vet{y})\in\reali^{n_2}\times\reali^{n_2}$.

Our aim is to determine a canonical transformation
$(\vet{p},\vet{q},\vet{x},\vet{y})=
\Kscr^{(\infty)}(\vet{P},\vet{Q},\vet{X},\vet{Y})$ which gives the
Hamiltonian $H^{(\infty)}=H^{(0)}\circ\Kscr^{(\infty)}$ the normal
form\footnote{We introduce a minor change of notations, using the
  symbol ${\omega}$ for the frequency vector related to the motion on
  a torus and ${\Omega}$ for the frequencies of the oscillations
  transverse to the elliptic torus.  The reader should avoid confusion
  with the previous use, when $\omega$ and $\Omega$ denoted the
  longitudes of the perihelia and of the nodes, respectively, which is
  the classical notation in Celestial Mechanics.}
\begin{equation}
\vcenter{\openup1\jot\halign{
 \hbox {\hfil $\displaystyle {#}$}
&\hbox {\hfil $\displaystyle {#}$\hfil}
&\hbox {$\displaystyle {#}$\hfil}\cr
H^{(\infty)}(\vet{P},\vet{Q},\vet{X},\vet{Y}) &=
&\vet{\omega}^{(\infty)}\cdot\vet{P}+
\sum_{j=1}^{n_2}\frac{\Omega_j^{(\infty)}\left(X_j^2+Y_j^2\right)}{2}+
\cr
& &\Oscr\big(\|\vet{P}\|^2\big)+
\Oscr\big(\|\vet{P}\|\|(\vet{X},\vet{Y})\|\big)+
\Oscr\big(\|(\vet{X},\vet{Y})\|^3\big)\,.
\cr
}}
\label{Ham^(infty)}
\end{equation}
Here the notation means that we want to remove all terms which are
linear in $\vet{P}$ and independent of $(\vet{X},\vet{Y})\,$, or at
most quadratic in $(\vet{X},\vet{Y})$ and independent of $\vet{P}$.  A
solution of the Hamilton's equations is
\begin{equation}
\big(\vet{P}(t),\vet{Q}(t),\vet{X}(t),\vet{Y}(t)\big)=
\big(\vet{0},\vet{Q}_0+\vet{\omega}^{(\infty)}t,\vet{0},\vet{0}\big)\,.
\label{sol_su_toro_ell}
\end{equation}
Choosing the initial conditions
$(\vet{P},\vet{Q},\vet{X},\vet{Y})=(\vet{0},\vet{Q}_0,\vet{0},\vet{0})$
(with $\vet{Q}_0\in\toro^{n_1}$), one immediately sees that the
corresponding orbit lies on the invariant $n_1$-dimensional torus
$\vet{P}=\vet{0}$ and $\vet{X}=\vet{Y}=\vet{0}\,$ and that the orbits
are quasi-periodic on it with frequencies $\vet{\omega}^{(\infty)}\,$.

The generic $r$-th step of our normalization algorithm is performed as
follows.  We reorder the Hamiltonian as
\begin{equation}
H^{(r-1)}(\vet{p},\vet{q},\vet{x},\vet{y})=
\vet{\omega}^{(r-1)}\cdot\vet{p}+\vet{\Omega}^{(r-1)}\cdot\vet{J}+
\sum_{s=0}^{\infty}\sum_{l=0}^{\infty}\sum_{{2j_1+j_2=l}\atop{j_1\ge 0\,,\,j_2\ge 0}}
f_{j_1\,,\,j_2}^{(r-1,s)}(\vet{p},\vet{q},\vet{x},\vet{y}) \,,
\label{Ham^(r-1)}
\end{equation}
where $J_j=(x_j^2+y_j^2)/2$ is the action which is usually related to
the $j$-th pair of secular canonical coordinates
$(x_j,y_j)\,,\ \forall\ j=1\,,\,\ldots\,,\,n_2\,$.  Moreover we pick a
suitable integer $K>0$ and determine the functions
$f_{j_1\,,\,j_2}^{(r-1,s)}$ so that they satisfy the following rules:
\begin{enumerate}
\renewcommand{\itemsep}{0pt}
\item[(A)] $f_{j_1\,,\,j_2}^{(r-1,s)}\in\Pscr_{j_1\,,\,j_2}^{(sK)}\,$,
  where $\Pscr_{j_1\,,\,j_2}^{(sK)}$ is the class of functions which
   are homogeneous polynomials of degree $j_1$ in the
  actions $\vet{p}\,$, homogeneous polynomials of degree
  $j_2$ in the secular variables $(\vet{x},\vet{y})\,$, and
  trigonometric polynomials of degree $sK$ in the angles $\vet{q}\,$;
\item[(B)] the terms $f_{j_1\,,\,j_2}^{(r-1,s)}$ are ``well
  Fourier-ordered''; this nonstandard definition means that
  $\forall\ j_1\ge 0\,,\,j_2\ge 0\,,\,s\ge 1\,$ every Fourier
  harmonic~$\vet{k}$ appearing in the expansion of
  $f_{j_1\,,\,j_2}^{(r-1,s)}$ is such that its corresponding
  trigonometric degree $|\vet{k}|=|k_1|+\ldots+|k_{n_1}|$ satisfies
  $(s-1)K<|\vet{k}|\leq sK\,$.
\end{enumerate}

By using formula~(\ref{Ham-trasl-fast}) and the properties~(i)--(iii)
of the canonical transformation $\Dscr\,$, one easily sees that the
Hamiltonian $H^{(0)}$ defined in~(\ref{Ham^(0)}) can be expanded in
the form written in~(\ref{Ham^(r-1)}), after having suitably reordered
its Fourier expansion so to satisfy the above requirements~(A)
and~(B). Therefore, our constructive algorithm can be applied to the
Hamiltonian $H^{(0)}$ by starting with $r=1\,$.

The comparison of the expansion in~(\ref{Ham^(r-1)}) with the normal
form in~(\ref{Ham^(infty)}) clearly shows that we should remove all
terms $f_{j_1\,,\,j_2}^{(0,s)}$ where the index $l=2j_1+j_2$ is such
that $0\le l\le 2\,$. Thus, the $r$-th step of our algorithm can be
naturally divided in three stages, each one aiming to reduce the
perturbation terms with $l=0,\,1,\,2$, respectively.

\subsection{First Stage of the  Normalization Step}\label{sbs:firststep}
 
The aim of the first stage in the $r$-th normalization step is to
remove the terms depending only on $\vet{q}\,$.

We use the Lie series algorithm to calculate the canonical
transformations (see, e.g.,~\cite{Giorgilli-1995} for an
introduction).  The generating function
$\chi_0^{(r)}(\vet{q})\in\Pscr_{0,0}^{(rK)}$ is determined by solving
the equation
\begin{equation}
\adaptpoisson{\chi_0^{(r)}}{\vet{\omega}^{(r-1)}\cdot\vet{p}}+
\sum_{s=1}^{r} f_{0,0}^{(r-1,s)}(\vet{q})
= 0\,,
\label{eq_per_chi_0}
\end{equation}
and the new Hamiltonian is determined as $H^{({\rm
    I};r)}=\exp\lie{\chi_0^{(r)}}H^{(r-1)}$, the symbol
$\poisson{\cdot}{\cdot}$ denoting the Poisson bracket.  The equation
for the generating function admits a solution provided the frequency
vector $\vet{\omega}^{(r-1)}$ is non-resonant up to order $rK\,$.
More precisely, we assume that
\begin{equation}
\min_{0<|\vet{k}|\le rK}\big|\vet{k}\cdot\vet{\omega}^{(r-1)}\big|\ge \alpha_{r}
\qquad{\rm with\ \ }\alpha_{r}>0\,,
\label{cond_non_ris_freq_vel}
\end{equation}
where $\{\alpha_{r}\}_{r>0}$ is a sequence of real positive numbers
and $|\vet{k}|=|k_1|+\ldots+|k_{n_1}|\,$.  The solution of the
homological equation~(\ref{eq_per_chi_0}) can be easily recovered by
looking at the little more complicate case of $X_2^{(r)}$, which is
discussed in the third stage of the $r$-th normalization step (see
formulas~(\ref{eq_per_X_2})--(\ref{exp_f_10})).

The canonical transformation writes $(\vet{p},\vet{q},\vet{x},\vet{y})
=\exp\lie{\chi_0^{(r)}}(\vet{p'},\vet{q'},\vet{x'},\vet{y'})\,$.
Omitting primes in order to simplify the notations the transformed
Hamiltonian can be written as
\begin{equation}
H^{({\rm I};r)}(\vet{p},\vet{q},\vet{x},\vet{y})=
\vet{\omega}^{(r-1)}\cdot\vet{p}+\vet{\Omega}^{(r-1)}\cdot\vet{J}+
\sum_{s=0}^{\infty}\sum_{l=0}^{\infty}\sum_{{2j_1+j_2=l}\atop{j_1\ge 0\,,\,j_2\ge 0}}
f_{j_1\,,\,j_2}^{({\rm I};r,s)}(\vet{p},\vet{q},\vet{x},\vet{y})\,.
\label{Ham^(I;r)}
\end{equation}
The functions $f_{j_1\,,\,j_2}^{({\rm I};r,s)}$ are recursively
defined, we omit this lengthy calculation since it is straightforward.
The main remark is concerned with the classes of functions.  It is
easy to check that
\begin{equation}
\frac{1}{i!}\lie{\chi_0^{(r)}}^i f_{j_1\,,\,j_2}^{(r-1,s)}
\in\Pscr_{j_1-i\,,\,j_2}^{((s+ir)K)}
\qquad
\forall\ 0\le i\le j_1\,,\ j_2\ge 0\,,\ s\ge 0\,,
\label{rel_classi_chi0}
\end{equation}
but it does not satisfy condition (B) at the beginning of the present
section.  Therefore, after having constructed $\sum_{j=0}^{s+ir}
f_{j_1-i\,,\,j_2}^{({\rm I};r,j)}$ by calculating all Poisson brackets
in $\frac{1}{i!}\lie{\chi_0^{(r)}}^i f_{j_1\,,\,j_2}^{(r-1,s)}$, we
perform a suitable reordering of the Taylor-Fourier series, so that
the expansion~(\ref{Ham^(I;r)}) satisfies both conditions~(A) and~(B).

\subsection{Second Stage of the  Normalization Step}\label{sbs:secondstep}

The aim of the second stage of the $r$-th normalization step is to
remove the terms which are linear in $(\vet{x},\vet{y})$ and
independent of $\vet{p}\,$.

We construct a second generating function
$\chi_1^{(r)}(\vet{q},\vet{x},\vet{y})\in\Pscr_{0,1}^{(rK)}$ by
solving the equation
\begin{equation}
\adaptpoisson{\chi_1^{(r)}}{\vet{\omega}^{(r-1)}\cdot\vet{p}+
\sum_{j=1}^{n_2}\frac{\Omega_j^{(r-1)}}{2}\big(x_j^2+y_j^2\big)}+
\sum_{s=0}^{r} f_{0,1}^{({\rm I};r,s)}(\vet{q},\vet{x},\vet{y})
= 0\,.
\label{eq_per_chi_1}
\end{equation}
With this generating function we construct the new Hamiltonian
$H^{({\rm II};r)}=\exp\lie{\chi_1^{(r)}}H^{({\rm I};r)}\,$.

Let us show how equation~(\ref{eq_per_chi_1}) is solved.  It is
convenient to temporarily introduce action-angle coordinates in place
of the secular pairs $(\vet{x},\vet{y})$ by putting
$x_j=\sqrt{2J_j}\cos\phi_j$ and $y_j=\sqrt{2J_j}\sin\phi_j$
$\forall\ j=1,\,\ldots\,,\,n_2\,$, so that the expansion of the known
terms appearing in equation~(\ref{eq_per_chi_1}) has the form:
\begin{equation}
\sum_{s=0}^{r} f_{0,1}^{({\rm I};r,s)}(\vet{q},\vet{J},\vet{\phi})
=\sum_{0\le|\vet{k}|\le rK}\sum_{j=1}^{n_2}
\sqrt{2J_j}\left[
c_{\vet{k}\,,\,j}^{(\pm)}\cos\big(\vet{k}\cdot\vet{q}\pm\phi_j\big)+
d_{\vet{k}\,,\,j}^{(\pm)}\sin\big(\vet{k}\cdot\vet{q}\pm\phi_j\big)\right]\,,
\label{exp_f_01}
\end{equation}
with known real coefficients $c_{\vet{k}\,,\,j}^{(\pm)}$ and
$d_{\vet{k}\,,\,j}^{(\pm)}\,$. Thus, one can easily check that
\begin{equation}
\chi_1^{(r)}(\vet{q},\vet{J},\vet{\phi})
=\sum_{0\le|\vet{k}|\le rK}\sum_{j=1}^{n_2}
\sqrt{2J_j}\left[
-\frac{c_{\vet{k}\,,\,j}^{(\pm)}\sin\big(\vet{k}\cdot\vet{q}\pm\phi_j\big)}
{\vet{k}\cdot\vet{\omega}^{(r-1)}\pm\Omega_j^{(r-1)}}+
\frac{d_{\vet{k}\,,\,j}^{(\pm)}\cos\big(\vet{k}\cdot\vet{q}\pm\phi_j\big)}
{\vet{k}\cdot\vet{\omega}^{(r-1)}\pm\Omega_j^{(r-1)}}\right]
\label{sol_per_chi_1}
\end{equation}
is a solution of the homological equation~(\ref{eq_per_chi_1}), and it
is consistently constructed provided the frequency vector
$\vet{\omega}^{(r-1)}$ satisfies the so-called {\it first Melnikov
  non-resonance condition\/} up to order $rK\,$, i.e.,
\begin{equation}
\min_{{\scriptstyle 0<|\vet{k}|\le rK}\atop{\scriptstyle j=1\,,\,\ldots\,,\,n_2}}
\big|\vet{k}\cdot\vet{\omega}^{(r-1)}\pm\Omega_j^{(r-1)}\big|
\ge \alpha_{r}
\qquad{\rm with\ \ }\alpha_{r}>0\,,
\label{cond_non_ris_Melnikov_1}
\end{equation}
and all the entries of the frequency vector $\vet{\Omega}^{(r-1)}$ are
far enough from the origin, i.e.,
\begin{equation}
\min_{j=1\,,\,\ldots\,,\,n_2}
\big|\Omega_j^{(r-1)}\big|
\ge \beta \qquad{\rm with\ \ }\beta>0\,.
\label{cond_non_ris_Omega_1}
\end{equation}
We remark that the latter condition~(\ref{cond_non_ris_Omega_1}) is
needed in general, but it is not necessary in the case of a planetary
Hamiltonian.  For, d'Alembert rules hold true, and so all
coefficients $c_{\vet{k}\,,\,j}^{(\pm)}$ and
$d_{\vet{k}\,,\,j}^{(\pm)}$ appearing in~(\ref{exp_f_01}) and having
even values of $|\vet{k}|$ are equal to zero.  However, such condition
is substantially included in a further one
(i.e.,~(\ref{cond_non_ris_Omega_2})) that we will need to introduce
later.

Starting from the expansion~(\ref{sol_per_chi_1}) of
$\chi_1^{(r)}(\vet{q},\vet{J},\vet{\phi})\,$, one can immediately
recover the expression of $\chi_1^{(r)}(\vet{q},\vet{x},\vet{y})$ as a
function of the original polynomial variables. We can then explicitly
calculate the expansion of the new Hamiltonian, which is written as
\begin{equation}
H^{({\rm II};r)}(\vet{p},\vet{q},\vet{x},\vet{y})=
\vet{\omega}^{(r-1)}\cdot\vet{p}+\vet{\Omega}^{(r-1)}\cdot\vet{J}+
\sum_{s=0}^{\infty}\sum_{l=0}^{\infty}\sum_{{2j_1+j_2=l}\atop{j_1\ge 0\,,\,j_2\ge 0}}
f_{j_1\,,\,j_2}^{({\rm II};r,s)}(\vet{p},\vet{q},\vet{x},\vet{y})\,.
\label{Ham^(II;r)}
\end{equation}
Here too we omit the explicit recursive expressions of the terms
$f_{j_1\,,\,j_2}^{({\rm II};r,s)}$, which can be obtained by a quite
annoying calculation.  We just add a remark on the implementation via
computer algebra.  Let us remark that
\begin{equation}
\frac{1}{i!}\lie{\chi_1^{(r)}}^i \sum_{2j_1+j_2=l}f_{j_1\,,\,j_2}^{({\rm I};r,s)}
\in\bigcup_{2j_1+j_2=l-i}\Pscr_{j_1\,,\,j_2}^{((s+ir)K)}
\qquad
\forall\ 0\le i\le l\,,\ s\ge 0\,.
\label{rel_classi_chi1}
\end{equation}
Therefore we construct the sum $\sum_{j=0}^{s+ir} \sum_{2j_1+j_2=l-i}
f_{j_1\,,\,j_2}^{({\rm II};r,j)}\,$ by calculating all the Poisson
brackets appearing in the expression of
$\frac{1}{i!}\lie{\chi_1^{(r)}}^i\sum_{2j_1+j_2=l}f_{j_1\,,\,j_2}^{({\rm
    I};r,s)}\,$ and the we proceed again with a reordering of the
Taylor-Fourier series so that also the expansion~(\ref{Ham^(II;r)})
satisfies the conditions~(A) and~(B), which have been stated at the
beginning of the present section.

\subsection{Third Stage of the  Normalization Step}\label{sbs:thirdstep}

The third and last stage of the $r$-th normalization step is more
elaborated.  The aim is to remove two classes of terms, namely terms
which are linear in $\vet{p}$ and independent of
$(\vet{x},\vet{y})\,$, and terms which are quadratic in
$(\vet{x},\vet{y})$ and independent of $\vet{p}\,$.

The Hamiltonian produced at the end of the $r$-th normalization step
is provided by the composition of three canonical
transformations\footnote{When one focuses on the estimates needed to
  prove the convergence of the algorithm, it is certainly simpler to
  introduce a single generating function
  $\chi_2^{(r)}(\vet{p},\vet{q},\vet{x},\vet{y})=
  X_2^{(r)}(\vet{p},\vet{q})+Y_2^{(r)}(\vet{q},\vet{x},\vet{y})$ and
  to consider the new Hamiltonian
  $\exp\lie{\Dscr_2^{(r)}}\circ\exp\lie{\chi_2^{(r)}}H^{({\rm II};r)}$
  which slightly differs from $H^{(r)}$, because $X_2^{(r)}$, and
  $Y_2^{(r)}$ do not commute with respect to the Poisson brackets.
  However, in the present work we prefer to split the present third
  stage of the $r$-th normalization step in three parts, so to
  highlight their different roles. Moreover, this choice looks more
  natural when one implements the constructive algorithm by algebraic
  manipulations on a computer.}  which can be given in terms of Lie
series, thus constructing the final Hamiltonian
$H^{(r)}=\exp\lie{\Dscr_2^{(r)}}\circ\exp\lie{Y_2^{(r)}}\circ
\exp\lie{X_2^{(r)}}H^{({\rm II};r)}\,$.  The generating functions
belong to three different classes:
$X_2^{(r)}(\vet{p},\vet{q})\in\Pscr_{1,0}^{(rK)}\,$,
$Y_2^{(r)}(\vet{q},\vet{x},\vet{y})\in\Pscr_{0,2}^{(rK)}$ and
$\Dscr_2^{(r)}(\vet{x},\vet{y})\in\Pscr_{0,2}^{(0)}\,$.  The explicit
expressions of these generating functions are given below, in
formulas~(\ref{sol_per_X_2}), (\ref{sol_per_Y_2})
and~(\ref{sol_per_D_2}), respectively.

We start with $X_2^{(r)}(\vet{p},\vet{q})\in\Pscr_{1,0}^{(rK)}\,$,
which is determined as the solution of the equation
\begin{equation}
\adaptpoisson{X_2^{(r)}}{\vet{\omega}^{(r-1)}\cdot\vet{p}}+
\sum_{s=1}^{r} f_{1,0}^{({\rm II};r,s)}(\vet{p},\vet{q})
= 0\,.
\label{eq_per_X_2}
\end{equation}
Writing the expansion of the known terms
appearing in equation~(\ref{eq_per_X_2}) in the form
\begin{equation}
\sum_{s=1}^{r} f_{1,0}^{({\rm II};r,s)}(\vet{p},\vet{q})
=\sum_{0<|\vet{k}|\le rK}\sum_{j=1}^{n_1}
p_j\left[
c_{\vet{k}\,,\,j}\cos\big(\vet{k}\cdot\vet{q}\big)+
d_{\vet{k}\,,\,j}\sin\big(\vet{k}\cdot\vet{q}\big)\right]\,,
\label{exp_f_10}
\end{equation}
with known real coefficients $c_{\vet{k}\,,\,j}$ and
$d_{\vet{k}\,,\,j}$ the generating function is given by
\begin{equation}
X_2^{(r)}(\vet{p},\vet{q})
=\sum_{0<|\vet{k}|\le rK}\sum_{j=1}^{n_1}
p_j\left[
-\frac{c_{\vet{k}\,,\,j}\sin\big(\vet{k}\cdot\vet{q}\big)}
{\vet{k}\cdot\vet{\omega}^{(r-1)}}+
\frac{d_{\vet{k}\,,\,j}\cos\big(\vet{k}\cdot\vet{q}\big)}
{\vet{k}\cdot\vet{\omega}^{(r-1)}}\right]\,,
\label{sol_per_X_2}
\end{equation}
The latter expression is consistent if the frequency vector
$\vet{\omega}^{(r-1)}$ satisfies the non-resonance
condition~(\ref{cond_non_ris_freq_vel}).

The generating function
$Y_2^{(r)}(\vet{q},\vet{x},\vet{y})\in\Pscr_{0,2}^{(rK)}$ is
determined by solving the equation
\begin{equation}
\adaptpoisson{Y_2^{(r)}}{\vet{\omega}^{(r-1)}\cdot\vet{p}+
\sum_{j=1}^{n_2}\frac{\Omega_j^{(r-1)}}{2}\big(x_j^2+y_j^2\big)}+
\sum_{s=1}^{r} f_{0,2}^{({\rm II};r,s)}(\vet{q},\vet{x},\vet{y})
= 0\,.
\label{eq_per_Y_2}
\end{equation}
Proceeding as in the second stage we use action-angle variables
$(\vet{J},\vet{\phi})$ in place of the secular pairs
$(\vet{x},\vet{y})\,$, so that the expansion of the known terms
appearing in equation~(\ref{eq_per_Y_2}) takes the form
\begin{equation}
\vcenter{\openup1\jot\halign{
 \hbox {\hfil $\displaystyle {#}$}
&\hbox {\hfil $\displaystyle {#}$\hfil}
&\hbox {$\displaystyle {#}$\hfil}\cr
\sum_{s=1}^{r} f_{0,2}^{({\rm II};r,s)}(\vet{q},\vet{J},\vet{\phi})
&=
&\sum_{0<|\vet{k}|\le rK}\,\sum_{i\,,\,j\,=\,1}^{n_2}\,
2\sqrt{J_iJ_j}\bigg[
c_{\vet{k}\,,\,i\,,\,j}^{(\pm\,,\,\pm)}
\cos\big(\vet{k}\cdot\vet{q}\pm\phi_i\pm\phi_j\big)
\cr
&
&\phantom{\sum_{0\le|\vet{k}|\le rK}\,\sum_{i\,,\,j\,=\,1}^{n_2}\,2\sqrt{J_iJ_j}\bigg[}
+d_{\vet{k}\,,\,i\,,\,j}^{(\pm\,,\,\pm)}
\sin\big(\vet{k}\cdot\vet{q}\pm\phi_i\pm\phi_j\big)\bigg]\,,
\cr
}}
\label{exp_f_02}
\end{equation}
with known real coefficients $c_{\vet{k}\,,\,i\,,\,j}^{(\pm\,,\,\pm)}$ and
$d_{\vet{k}\,,\,i\,,\,j}^{(\pm\,,\,\pm)}\,$.  Hence the generating
function is determined as 
\begin{equation}
\vcenter{\openup1\jot\halign{
 \hbox {\hfil $\displaystyle {#}$}
&\hbox {\hfil $\displaystyle {#}$\hfil}
&\hbox {$\displaystyle {#}$\hfil}\cr
Y_2^{(r)}(\vet{q},\vet{J},\vet{\phi})
&=
&\sum_{0<|\vet{k}|\le rK}\,\sum_{i\,,\,j\,=\,1}^{n_2}\,
2\sqrt{J_iJ_j}\Bigg[
-\frac{c_{\vet{k}\,,\,i\,,\,j}^{(\pm\,,\,\pm)}
\sin\big(\vet{k}\cdot\vet{q}\pm\phi_i\pm\phi_j\big)}
{\vet{k}\cdot\vet{\omega}^{(r-1)}\pm\Omega_i^{(r-1)}\pm\Omega_j^{(r-1)}}
\cr
&
&\phantom{\sum_{0<|\vet{k}|\le rK}\,\sum_{i\,,\,j\,=\,1}^{n_2}\,2\sqrt{J_iJ_j}\bigg[}
+\frac{d_{\vet{k}\,,\,i\,,\,j}^{(\pm\,,\,\pm)}
\cos\big(\vet{k}\cdot\vet{q}\pm\phi_i\pm\phi_j\big)}
{\vet{k}\cdot\vet{\omega}^{(r-1)}\pm\Omega_i^{(r-1)}\pm\Omega_j^{(r-1)}}
\Bigg]\,.
\cr
}}
\label{sol_per_Y_2}
\end{equation}
In order to make this expression consistent we must assume that the
frequency vector $\vet{\omega}^{(r-1)}$ satisfies the so-called {\it
  second Melnikov non-resonance condition\/} up to order $rK\,$, i.e.,
\begin{equation}
\min_{{\scriptstyle 0<|\vet{k}|\le rK}\atop{\scriptstyle
    i\,,\,j\,=1\,,\,\ldots\,,\,n_2}}
\big|\vet{k}\cdot\vet{\omega}^{(r-1)}\pm\Omega_i^{(r-1)}\pm\Omega_j^{(r-1)}\big|
\ge \alpha_{r} \qquad{\rm with\ \ }\alpha_{r}>0\,.
\label{cond_non_ris_Melnikov_2}
\end{equation}
Let us here remark that the latter condition includes also the
non-resonance condition~(\ref{cond_non_ris_freq_vel}) as a special
case, i.e., when $i=j$ and the signs appearing in the expression
$\pm\Omega_i^{(r-1)}\pm\Omega_j^{(r-1)}$ are opposite.

Concerning the the generating function $\Dscr_2^{(r)}$, once again it
is convenient to replace the secular pairs $(\vet{x},\vet{y})$ with
the action-angle coordinates $(\vet{J},\vet{\phi})\,$.  Remark that
$\vet{\Omega}^{(r-1)}\cdot\vet{J}$ and $f_{0\,,\,2}^{({\rm
    II};r,0)}(\vet{x},\vet{y})$ are the only terms appearing in
expansion~(\ref{Ham^(II;r)}) which are quadratic in
$(\vet{x},\vet{y})$ (so they also are~$\Oscr(\vet{J})\,$) and do not
depend on $\vet{p}$ and $\vet{q}\,$.  The canonical transformation
induced by the Lie series $\exp\lie{\Dscr_2^{(r)}}$ aims to eliminate
the part of $f_{0\,,\,2}^{({\rm II};r,0)}$ depending on the secular
angles $\vet{\phi}\,$.  Therefore, the generating function
$\Dscr_2^{(r)}$ is determined as the solution of the equation
\begin{equation}
\adaptpoisson{\Dscr_2^{(r)}}{\vet{\Omega}^{(r-1)}\cdot\vet{J}}+
f_{0\,,\,2}^{({\rm II};r,0)}(\vet{J},\vet{\phi})-
\big\langle f_{0\,,\,2}^{({\rm II};r,0)}\big\rangle_{\vet{\phi}} = 0\,,
\label{eq_per_D_2}
\end{equation}
where $\langle\cdot\rangle_{\vet{\phi}}$ denotes the average with
respect to the angles $\vet{\phi}\,$.  Writing the expansion of the
known terms as
\begin{equation}
\vcenter{\openup1\jot\halign{
 \hbox {\hfil $\displaystyle {#}$}
&\hbox {\hfil $\displaystyle {#}$\hfil}
&\hbox {$\displaystyle {#}$\hfil}\cr
f_{0,2}^{({\rm II};r,0)}(\vet{J},\vet{\phi})
&=
&\sum_{i\,,\,j\,=\,1}^{n_2}\,
\sum_{{{\scriptstyle s_i=\pm 1}}\atop{{\scriptstyle s_j=\pm 1}}}\,
2\sqrt{J_iJ_j}\bigg[
c_{i\,,\,j\,,\,s_i\,,\,s_j}
\cos\big(s_i\phi_i+s_j\phi_j\big)
\cr
& &\phantom{\sum_{i\,,\,j\,=\,1}^{n_2}\,
\sum_{{{\scriptstyle s_i=\pm 1}}\atop{{\scriptstyle s_j=\pm 1}}}\,
2\sqrt{J_iJ_j}\bigg[}
+d_{i\,,\,j\,,\,s_i\,,\,s_j}
\sin\big(s_i\phi_i+s_j\phi_j\big)\bigg]\,,
\cr
}}
\label{exp_f_02_mediata}
\end{equation}
with known real coefficients $c_{i\,,\,j\,,\,s_i\,,\,s_j}$ and
$d_{i\,,\,j\,,\,s_i\,,\,s_j}$ the generating function is given by
\begin{equation}
\vcenter{\openup1\jot\halign{
 \hbox {\hfil $\displaystyle {#}$}
&\hbox {\hfil $\displaystyle {#}$\hfil}
&\hbox {$\displaystyle {#}$\hfil}\cr
\Dscr_{2}^{(r)}(\vet{J},\vet{\phi})
&= &\sum_{i\,,\,j\,=\,1}^{n_2}\,
\sum_{{{\scriptstyle s_i\,,\,s_j\,=\,\pm 1}}
\atop{{\scriptstyle i\cdot s_i + j\cdot s_j\,\neq\,0}}}\,2\sqrt{J_iJ_j}\Bigg[
-\frac{c_{i\,,\,j\,,\,s_i\,,\,s_j}\sin\big(s_i\phi_i+s_j\phi_j\big)}
{s_i\Omega_i^{(r-1)}+s_j\Omega_j^{(r-1)}}
\cr
& &\phantom{\sum_{i\,,\,j\,=\,1}^{n_2}\,
\sum_{{{\scriptstyle s_i\,,\,s_j\,=\,\pm 1}}
\atop{{\scriptstyle s_i\cdot i + s_j\cdot j\,\neq\,0}}}\,2\sqrt{J_iJ_j}\Bigg[}
+\frac{d_{i\,,\,j\,,\,s_i\,,\,s_j}\cos\big(s_i\phi_i+s_j\phi_j\big)}
{s_i\Omega_i^{(r-1)}+s_j\Omega_j^{(r-1)}}
\Bigg]\,.
\cr
}} 
\label{sol_per_D_2}
\end{equation}
The solution is consistent provided the frequency vector
$\vet{\Omega}^{(r-1)}$ satisfies the finite non-resonance condition
\begin{equation}
\min_{|\vet{l}|=2}
\big|\vet{l}\cdot\vet{\Omega}^{(r-1)}\big|
\ge \beta \qquad{\rm with\ \ }\beta>0\,.
\label{cond_non_ris_Omega_2}
\end{equation}
Having determined the generating functions $Y_2^{(r)}$ and
$\Dscr_{2}^{(r)}$ it is a standard matter to replace the action-angle
variables $\vet{J},\vet{\phi}$ with the secular variables
$\vet{x},\vet{y}\,$.

At this point of the algorithm, it is convenient to slightly modify
the frequencies $\vet{\omega}^{(r-1)}$ and $\vet{\Omega}^{(r-1)}$, so
to include the terms which are linear with respect to the actions and
do not depend on the angles.  Such terms can not be eliminated by our
normalization procedure. More precisely, we define
$\vet{\omega}^{(r)}$ and $\vet{\Omega}^{(r)}\,$, so that
\begin{equation}
\vet{\omega}^{(r)}\cdot\vet{p}=
\vet{\omega}^{(r-1)}\cdot\vet{p}+
f_{1,0}^{({\rm II};r,0)}(\vet{p})\,,
\qquad
\vet{\Omega}^{(r)}\cdot\vet{J}=
\vet{\Omega}^{(r-1)}\cdot\vet{J}+
\big\langle f_{0\,,\,2}^{({\rm II};r,0)}\big\rangle_{\vet{\phi}}\,\,.
\label{aggiusta_frequenze}
\end{equation}

We are now able to explicitly produce the expansion of the new
Hamiltonian, which can be written as
\begin{equation}
H^{(r)}(\vet{p},\vet{q},\vet{x},\vet{y})=
\vet{\omega}^{(r)}\cdot\vet{p}+\vet{\Omega}^{(r)}\cdot\vet{J}+
\sum_{s=0}^{\infty}\sum_{l=0}^{\infty}\sum_{{2j_1+j_2=l}\atop{j_1\ge 0\,,\,j_2\ge 0}}
f_{j_1\,,\,j_2}^{(r,s)}(\vet{p},\vet{q},\vet{x},\vet{y})\,.
\label{Ham^(r)}
\end{equation}
Let us remark that this expansion of $H^{(r)}$ has exactly the same
form of that written for $H^{(r-1)}$ in~(\ref{Ham^(r-1)}), but we
stress that the algorithm is arranged so to make smaller and smaller
the contribution of the terms $f_{j_1\,,\,j_2}^{(r,s)}\,$, when the
value of $r$ is increased, $\forall\ s\ge 0$ and $l=2
j_1+j_2=0,\,1,\,2\,$.

Let us add a remark which may be useful in implementing the
transformation via computer algebra.  For what concerns the generating
function $X_2^{(r)}$ we have
\begin{equation}
\frac{1}{i!}\lie{X_2^{(r)}}^i f_{j_1\,,\,j_2}^{({\rm II};r,s)}
\in\Pscr_{j_1\,,\,j_2}^{((s+ir)K)}
\qquad
\forall\ i\ge 0\,,\ j_1\ge 0\,,\ \ j_2\ge 0\,,\ s\ge 0\,.
\label{rel_classi_X2}
\end{equation}
For the generating function $Y_2^{(r)}$ we have the relations
\begin{equation}
\frac{1}{i!}\lie{Y_2^{(r)}}^i \sum_{2j_1+j_2=l}f_{j_1\,,\,j_2}^{({\rm II};r,s)}
\in\bigcup_{2j_1+j_2=l}\Pscr_{j_1\,,\,j_2}^{((s+ir)K)}
\qquad
\forall\ i\ge 0\,,\ l\ge 0\,,\ s\ge 0\,.
\label{rel_classi_Y2}
\end{equation}
Finally, one can easily remark that each class of function is
invariant with respect to a Poisson bracket with the generating
function $\Dscr_2^{(r)}\,$.  Hence we have
\begin{equation}
\frac{1}{i!}\lie{\Dscr_2^{(r)}}^i f_{j_1\,,\,j_2}^{({\rm II};r,s)}
\in\Pscr_{j_1\,,\,j_2}^{(sK)}
\qquad
\forall\ i\ge 0\,,\ j_1\ge 0\,,\ \ j_2\ge 0\,,\ s\ge 0\,.
\label{rel_classi_D2}
\end{equation}
By taking into account the relations
(\ref{rel_classi_X2})--(\ref{rel_classi_D2}) among the classes of
functions, the definition~(\ref{aggiusta_frequenze}) of the new
frequencies vectors and by reordering the Taylor-Fourier series, it is
possible to ensure that also the expansion~(\ref{Ham^(r)}) satisfies
the conditions~(A) and~(B) which have been stated at the beginning of
the present section.

Thus the $r$-th step of the normalization procedure that we have
described here can be iterated.

\subsection{Some Remarks about the Convergence of the Normalization
  Algorithm}\label{sbs:remarks_theory}

We devote this section to an informal discussion of the relations
between the normalization procedure for an elliptic torus, which is
the subject of the present paper, and the Kolmogorov's algorithm for a
torus of maximal dimension.  Our aim is to bring into evidence, on the
one hand, the differences that make the case of an elliptic lower
dimensional torus definitely more difficult and, on the other hand,
the impact that these differences have on the explicit calculation.

The main hypotheses of Kolmogorov's theorem are (a)~that the
perturbation should be small enough and (b)~that a strong
non-resonance condition must be satisfied by the frequencies of the
unperturbed torus.  Both these conditions appear also in the proof of
existence of elliptic tori, but the condition of non-resonance
presents some critical peculiarities.

A common problem in perturbation theory is concerned with the
smallness of the perturbation, because the main analytical estimates
are usually extremely restrictive.  Nevertheless, we have strong
evidence that realistic estimates may be obtained by using algebraic
manipulations in order to implement a computer-assisted proof (see,
e.g.,~\cite{Loc-Gio-2000}).  In the case of Kolmogorov's theorem a
computer-assisted procedure takes advantage of the preliminary
application of the algorithm constructing the normal form (which is
explicitly performed for a {\it finite\/} number of steps $R\,$, as
large as possible), because a suitable version of the KAM theorem is
finally applied to the Hamiltonian $H^{(R)}$ having the perturbing
terms strongly reduced with respect to the initial $H^{(0)}$. In the
case of elliptic lower dimensional tori a similar procedure applies.
The explicit application of the normalization algorithm mainly
requires to translate into a programming language the method described
in the previous sections.  Concerning the actual reduction of the
perturbation, by comparing the Hamiltonian normal
form~(\ref{Ham^(infty)}) with the expansion~(\ref{Ham^(r)}) of
$H^{(r)}$, one easily realizes that the initial expression of the
perturbation (making part of the Hamiltonian $H^{(0)}$, written
in~(\ref{Ham^(0)})) is given by
\begin{equation}
\sum_{s=0}^{\infty}\sum_{l=0}^{2}\sum_{{2j_1+j_2=l}\atop{j_1\ge
    0\,,\,j_2\ge 0}} f_{j_1\,,\,j_2}^{(0,s)}\,.
\label{initial-pert-terms}
\end{equation}
Looking at all the preliminary expansions, which have been described
in sect.~\ref{sec:2D_plan_Ham} and allowed us to introduce the initial
Hamiltonian $H^{(0)}$, one immediately sees that all the perturbing
terms appearing in~(\ref{initial-pert-terms}) are proportional to
$\mu\,$. Let us also recall that the small parameter $\mu$ is equal to
the mass ratio between the biggest planet and the central star
(according to its definition given in the discussion following
formula~(\ref{Ham-iniz-Poincare-var})).

The non-resonance conditions on the frequencies represent a definitely
more delicate problem.  We recall that in order to apply Kolmogorov's
theorem for a maximal dimensional torus, one must choose the $n$
frequencies $\omega_1,\ldots,\omega_n$ so as to satisfy a strong
non-resonance condition.  A typical request is that they obey a
Diophantine condition, i.e., that the sequence
$\left\{\alpha_r\right\}_{r\ge 1}$ appearing in the
inequality~(\ref{cond_non_ris_freq_vel}) must be such that
$\alpha_r\ge\gamma/(rK)^\tau$ with suitable positive values of the
constant $\gamma$ and $\tau\,$.  This choice must be made at the very
beginning of the procedure, and the perturbed invariant torus that is
found at the end has {\it the same\/} frequencies as the unperturbed
one.  This is assured by introducing, at every step, a small translation
of the actions which keeps the frequencies constant.

In the case of the elliptic lower dimensional torus one deals instead
with two separate set of frequencies, namely
$\vet{\omega}^{(0)}\in\reali^{n_1}$ which characterize the orbits on
the torus, and the transverse frequencies
$\vet{\Omega}^{(0)}\in\reali^{n_2}$ that are related to the
oscillations of orbits close to but not lying on the torus.  Now, the
frequencies $\vet{\omega}^{(0)}$ on the torus can be chosen in an
arbitrary way, but the transverse frequencies $\vet{\Omega}^{(0)}$ are
functions of $\vet{\omega}^{(0)}$, being given by the Hamiltonian.
This is easily understood by considering the case of a periodic orbit,
i.e., $n_1=1$, since in that case the transverse frequencies are
related to the eigenvalues of the monodromy matrix.

The striking fact is that, due precisely to the dependence of the
transverse frequencies $\vet{\Omega}^{(0)}$ on $\vet{\omega}^{(0)}$,
the algorithm forces us to change these frequencies at every step.
That is, one actually deals with infinite sequences
$\vet{\omega}^{(r)}$ and $\vet{\Omega}^{(r)}$, all required to satisfy
at every order a non-resonance condition of the
form~(\ref{cond_non_ris_Melnikov_2}).  Moreover, both sequences should
converge to a final set of frequencies
$\vet{\omega}^{(\infty)}=\vet{\omega}^{(\infty)}\big(\vet{\omega}^{(0)}\big)$
and
$\vet{\Omega}^{(\infty)}=\vet{\Omega}^{(\infty)}\big(\vet{\omega}^{(0)}\big)$
which must be non-resonant (e.g., Diophantine).  Thus, we are forced
to conclude that, depending on the initial choice of
$\vet{\omega}^{(0)}$, it may happen that the algorithm stops at some
step because the frequencies fail to satisfy some of the non-resonance
conditions~(\ref{cond_non_ris_freq_vel}),
(\ref{cond_non_ris_Melnikov_1}), (\ref{cond_non_ris_Omega_1}),
(\ref{cond_non_ris_Melnikov_2}) and~(\ref{cond_non_ris_Omega_2}).
This is indeed one of the main difficulties in working out the proof
of existence of an elliptic torus.

Let us first give an informal description of the analytical procedure.
We leave the full exposition with the convergence proof to a
forthcoming paper.  One initially considers an open ball
$\Bscr\subset\reali^{n_1}$ such that the Diophantine condition {\it at
  finite order\/} required for the first step is satisfied by every
$\vet{\omega}^{(0)}\in\Bscr\subset\reali^{n_1}$ and by the
corresponding transverse frequencies $\vet{\Omega}^{(0)}$.  This can
be done, because only a finite number of non-resonance relations are
considered.  One must show that at every step there exists a subset of
frequencies in $\Bscr$ which satisfies the non-resonance conditions
(still at finite but increasing order) required in order to perform
the next step, together with the corresponding transverse frequencies.
This is obtained by a procedure which is inspired by Arnold's proof
scheme of KAM theorem, and is quite different from Kolmogorov's one
(compare~\cite{Kolmogorov-1954} with~\cite{Arnold-1963.1}).  At every
step one removes from $\Bscr$ a finite number of intersections of
$\Bscr$ with a small strip around a resonant plane in $\reali^{n_1}$,
assuring that the width of the strip decreases fast enough so that the
remaining set always has a non-empty interior part.  By the way, this
is also strongly reminiscent of the process of construction of a
Cantor set.  The final goal is to prove that one is left with a Cantor
set on non-resonant frequencies which satisfy the required resonance
conditions and has positive Lebesgue measure.  Moreover, the relative
measure with respect to $\Bscr$ tends to one when the size of the
perturbation is decreased to zero.  This is the underlying idea of the
proof that we plan to present in complete form
in~\cite{Gio-Loc-San-2011}.

Leaving again the proof of the algorithm convergence to a forthcoming
paper, we report more precise statements about the
previous informal discussion.  Let us emphasize that
Hamiltonian~\eqref{Ham^(r-1)} satisfies the following parity
condition: all coefficients of the expansion of $f_{l,m}^{(r,s)}$
having even (odd) trigonometrical degree in the $\vet{q}$ variables
and odd (even) degree in ($\vet{x}$, $\vet{y}$), are identically zero.
In the planetary case this is a straightforward consequence of
d'Alembert rules.

\begin{theorem}
Consider a set of Hamiltonians parametrized with respect to the
frequency $\vet{\omega}^{(0)}$ that are analytic on some suitable open
domain (of their canonical variables) and are of the form
\eqref{Ham^(r-1)} with $r=1\,$, namely
$$
H^{(0)}(\vet{p},\vet{q},\vet{x},\vet{y};\vet{\omega}^{(0)})=
\vet{\omega}^{(0)}\cdot\vet{p}+\vet{\Omega}^{(0)}(\vet{\omega}^{(0)})\cdot\vet{J}+
\sum_{s=0}^{\infty}\sum_{l=0}^{\infty}\sum_{{2j_1+j_2=l}\atop{j_1\ge 0\,,\,j_2\ge 0}}
f_{j_1\,,\,j_2}^{(0,s)}(\vet{p},\vet{q},\vet{x},\vet{y};\vet{\omega}^{(0)})\ .
$$
Assume the following hypotheses to be fulfilled:

\item[(i)] let
  $\vet{\omega}^{(0)}\mapsto\vet{\Omega}^{(0)}(\vet{\omega}^{(0)})$ be
  a regular function defined on an open set $\Bscr\subset\reali^{n_1}$
  such that
  $$
  |\vet{l}\cdot\vet{\Omega}^{(0)}(\vet{\omega}^{(0)})|\geq\Xi_0\ ,
  \ \ \forall\ |\vet{l}|=1,\,2\,,\ \vet{\omega}^{(0)}\in\Bscr\ ,
  $$
  and
  $$
  \|\vet{\Omega}^{(0)}(\vet{\omega}^{(0)})-
  \vet\Omega^{(0)}(\tilde{\vet\omega}^{(0)})\|\leq
  L_0\|\vet{\omega}^{(0)}-\tilde{\vet{\omega}}^{(0)}\|\ ,
  \ \ \forall\ \vet{\omega}^{(0)},\,\tilde{\vet{\omega}}^{(0)}\in\Bscr\ ,
  $$
  for some positive constants $\Xi_0$ and $L_0\,$;

\item[(ii)] the functions
  $f_{j_1\,,\,j_2}^{(0,s)}\in\Pscr_{j_1\,,\,j_2}^{(sK)}$
  $\forall\ j_1,\,j_2,\,s\geq0$ with a fixed integer value of $K>0\,$;
  furthermore, they satisfy the parity condition above and
  \begin{equation}
    \vcenter{\openup1\jot\halign{
        \hbox {\hfil $\displaystyle {#}$}
        &\hbox {\hfil $\displaystyle {#}$\hfil}
        &\hbox {$\displaystyle {#}$\hfil}\cr
    \qquad\qquad\ 
    f_{0,0}^{(0,0)}=f_{0,1}^{(0,0)}=f_{0,2}^{(0,0)}=f_{1,0}^{(0,0)} &= &0\ ,
    \cr
    \langle f_{0,0}^{(0,1)}\rangle_\vet{q} =\langle
    f_{0,1}^{(0,1)}\rangle_\vet{q} =\langle f_{0,2}^{(0,1)}\rangle_\vet{q}
    =\langle f_{1,0}^{(0,1)}\rangle_\vet{q} &= &0\ ;
    \cr
    }}
    \label{eq:parti_Ham_nulle}
  \end{equation}

\item[(iii)] the non-degeneracy condition
  $$
  \det\left(\frac{\partial^2f_{2\,,\,0}^{(0,0)}}{\partial p_i\partial p_j}
  \right)_{i,j=1,\ldots,n_1}\neq 0\ ,
  $$
  is satisfied $\forall\ \vet{\omega}^{(0)}\in\Bscr$ and $\vet{p}$ in the
  definition domain of the Hamiltonians;

\item[(iv)] the following upper bounds hold true:
  \begin{equation}
    \|f_{j_1\,,\,j_2}^{(0,s)}\|\leq \epsilon^s E\ ,
    \qquad\forall\ j_1,\,j_2,\,s\geq0\ ,
    \label{ineq_smallness-cond}
  \end{equation}
  for some positive constants $\epsilon$ and $E\,$, with $\epsilon$
  {\rm small enough}.

Then, there exists a Cantor set $\Sscr\subset\Bscr$ of positive
Lebesgue measure (in $\reali^{n_1}$), such that for each
$\vet{\omega}^{(0)}\in\Sscr$ there is an analytic canonical
transformation giving the Hamiltonian the normal form
\eqref{sol_su_toro_ell}, namely,
\begin{equation*}
\vcenter{\openup1\jot\halign{
 \hbox {\hfil $\displaystyle {#}$}
&\hbox {\hfil $\displaystyle {#}$\hfil}
&\hbox {$\displaystyle {#}$\hfil}\cr
H^{(\infty)}(\vet{P},\vet{Q},\vet{X},\vet{Y};\vet{\omega}^{(0)}) &=
&\vet{\omega}^{(\infty)}\cdot\vet{P}+
\sum_{j=1}^{n_2}\frac{\Omega_j^{(\infty)}\left(X_j^2+Y_j^2\right)}{2}+
\cr
& &\Oscr\big(\|\vet{P}\|^2\big)+
\Oscr\big(\|\vet{P}\|\|(\vet{X},\vet{Y})\|\big)+
\Oscr\big(\|(\vet{X},\vet{Y})\|^3\big)\ ,
\cr
}}
\end{equation*}
where $\vet{\omega}^{(\infty)}=\vet{\omega}^{(\infty)}(\vet{\omega}^{(0)})$
and $\vet{\Omega}^{(\infty)}=\vet{\Omega}^{(\infty)}(\vet{\omega}^{(0)})\,$.
\end{theorem}
In the statement above, the functional norm appearing
in~\eqref{ineq_smallness-cond} is defined as the sum of the absolute
values of the coefficients appearing in the Taylor-Fourier expansion
of $f_{j_1\,,\,j_2}^{(0,s)}\,$.  Let us recall that, as remarked at
the beginning of section 3, the $n_1$-dimensional (elliptic) torus
corresponding to $\vet{P}=\vet{0}$ and $\vet{X}=\vet{Y}=\vet{0}$ is
invariant, and carries quasi periodic orbits with frequencies
$\vet\omega^{(\infty)}$.  Moreover, let us emphasize that taking
control of the measure of the set $\Sscr$ is perhaps the main
difficulty of the proof. Indeed, the definition of $\Sscr$ is very
delicate, because it is such that
$\forall\ \vet{\omega}^{(0)}\in\Sscr$ the sequence
$\big\{\big(\vet{\omega}^{(r)},\,\vet{\Omega}^{(r)}\big)\big\}_{r\ge
  0}$ (that is introduced by the algorithm described in the present
section) starting from
$\big(\vet{\omega}^{(0)},\,\vet{\Omega}^{(0)}(\vet{\omega}^{(0)})\big)$
satisfies the following conditions $\forall\ r\ge 1\,$:
\begin{equation}
\min_{\scriptstyle{0<|\vet{k}|\leq rK}\atop\scriptstyle{0\leq|\vet{l}|\leq2}}|
\scalprod{\vet{k}}{\vet{\omega}^{(r-1)}}
+ \scalprod{\vet{l}}{\vet{\Omega}^{(r-1)}}|\geq\frac{\gamma}{(rK)^{\tau}}\ ,
\qquad
\min_{|\vet{l}|=2}|\scalprod{\vet{l}}{\vet{\Omega}^{(r-1)}}|\geq\beta
\label{eq_nnris}
\end{equation}
for some positive fixed values of $\gamma\,$, $\beta$
and $\tau > n_1-1\,$.

Let us now come to the numerical aspect.  At first sight the formal
algorithm seems to require a cumbersome trial and error procedure in
order to find the good frequencies: when some of the non-resonance
conditions fail to be satisfied at a given step, one should change the
initial frequencies and restart the whole process.  Moreover, since
the non-resonance condition must be satisfied by the final
frequencies, which obviously can not be calculated, the whole process
seems to be unsuitable for a rigorous proof.  We explain here in which
sense the computer-assisted proofs can help to improve the results
also in this context.  We make two remarks.

The first remark is connected with the use of interval arithmetic
while performing the actual construction.  Following the suggestion of
the analytic scheme of the proof, we look for uniform estimates on a
small open ball $\Bscr\,$, such that
$\forall\ \vet{\omega}^{(0)}\in\Bscr$ we explicitly perform $R$
normalization steps, with $R$ as large as possible.  Essentially, we
may reproduce numerically the process of eliminating step-by-step the
unwanted resonant frequencies by suitably determining the intervals.
Once $R$ steps have been explicitly performed, we may apply to the
partially normalized Hamiltonian $H^{(R)}$ a suitable formulation of
the KAM theorem for elliptic tori.  This means that we recover the
scheme that we have already applied to the case of maximal dimension
tori.  That is, we can take advantage of the fact that the perturbing
terms are much smaller than the corresponding ones for the initial
Hamiltonian $H^{(0)}$; thus, in principle we could ensure that for
realistic values of $\mu$ the relative measure of the invariant tori
is so large that the set of those $\vet{\omega}^{(0)}$ for which the
algorithm can not work (i.e., $\Bscr\setminus\Sscr$) is so small that
it can be neglected when we are dealing with a practical application.

The second remark is that we may take a more practical attitude,
relying on the fact that the set of good frequencies, according to the
theory, has Lebesgue measure close to one, so that the case of
frequencies which are resonant at some finite order occurs with very
low probability.  Thus, we just make a choice of the initial
frequencies and proceed with the construction, checking at every order
that the non-resonance conditions that we need at that order are
fulfilled.  We emphasize that the most extended resonant regions are
those of low order, so that it is not very difficult to check
initially that the chosen frequencies will likely be good enough.  It
may happen, of course, that the whole procedure must be restarted with
different frequencies, but this case is reasonably expected to occur
only rarely.  For, since the size of the perturbation is expected to
decrease geometrically, we may confidently expect that the probability
of failure will decrease, too.  This is actually confirmed by the
calculations we have done.

When $R$ steps have been made, in principle we can apply the theorem
to a small neighborhood of the calculated frequencies by choosing a
suitable initial ball around the frequencies approximated at that
step.

\section{Elliptic Tori for the SJSU System}\label{sec:num_tests}

We come now to the application of the formal algorithm for the
construction of an elliptic torus to the planar SJSU system.  Here we
explicitly construct the normal form at a finite order checking that
the norms of the generating functions decrease as predicted by the
theory.  Then we perform a numerical check by comparing the orbit
obtained via the normal form with the numerically integrated one.

The initial Hamiltonian is written in~(\ref{Ham-trasl-fast}), with a
suitable rearrangement of terms so that it is given the
form~(\ref{Ham^(r-1)}) with $r=1\,$.  This requires also a
diagonalization of the quadratic part in the secular variables, which
is performed as described at the end of sect.~\ref{sec:2D_plan_Ham}.

\subsection{Constructing the Elliptic Torus by Using Computer Algebra}
\label{sbs:comp_algebra}

\begin{figure}
\centerline{\includegraphics[width=140mm]{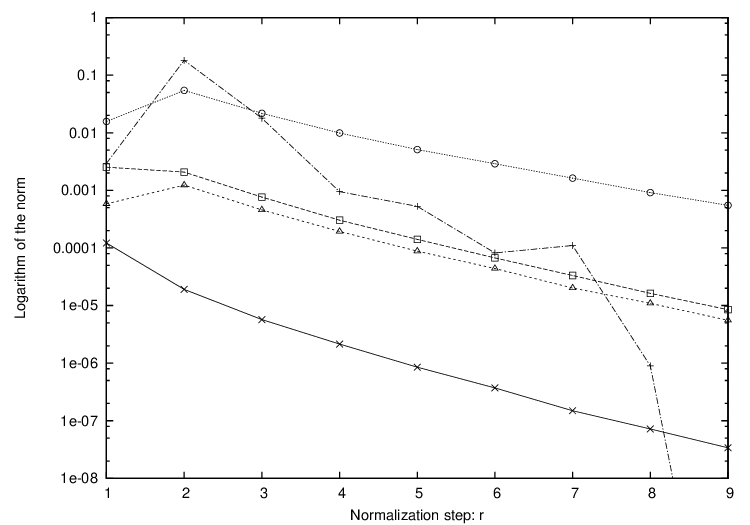}}
\caption{Numerical check of the algorithm constructing the normal form
  related to an elliptic torus for the planar SJSU system: plot of the
  norm of the generating functions as a function of the normalization
  step $r\,$; more precisely, the symbols $\times\,$, $\square\,$,
  $\triangle\,$, $\bigcirc$ and $+$ refer to the norm of the
  generating functions $\chi_0^{(r)}$, $\chi_1^{(r)}$, $X_2^{(r)}$,
  $Y_2^{(r)}$ and $\Dscr_2^{(r)}$, respectively, which are defined
  during the normalization algorithm, as described in
  sect.~\ref{sec:form_alg}. The norm is calculated by adding up the
  absolute values of all coefficients appearing in the expansion of
  each generating function.}
\label{fig:norm_gen_fun}
\end{figure}

We applied the algorithm constructing elliptic tori (which has been
widely described in sect.~\ref{sec:form_alg}) to the Hamiltonian
$H^{(0)}$ (that is defined in~(\ref{Ham^(0)}) and has been obtained as
described in sect.~\ref{sec:2D_plan_Ham}).  The parameters have been
fixed according to the specific values of the planar SJSU system,
which are reported in table~\ref{tab:parameters_2D_SJSU}.  Our
software package for computer algebra allowed us to explicitly
calculate all the expansions~(\ref{Ham^(0)}) of $H^{(r)}$ with index
$r$ ranging between~$0$ and $9\,$, so to include: (c1)~the terms
having degree $j_1$ in the actions $\vet{p}$ with $j_1\le 3\,$,
(c2)~all terms having degree~$j_2$ in the variables
$(\vet{x},\vet{y})\,$, with $j_2$ such that $2j_1+j_2\le 8\,$,
(c3)~all terms up to the trigonometric degree $18$ with respect to the
angles $\vet{q}\,$.  The truncation rules~(c1)--(c3) are in agreement
with those prescribed about the expansion~(\ref{Ham-trasl-fast}) in
sect.~\ref{sec:2D_plan_Ham} at points~(b1)--(b3).  Let us remark that
both the truncation rules~(c1) and~(c2) are preserved by all
canonical transformations included in our algorithm. Moreover, we have
found that fixing $K=2$ is a suitable choice to have a rather regular
decreasing of the size of the generating functions when the
normalization step $r$ is increased, as shown in
fig.~\ref{fig:norm_gen_fun}. Since the maximal trigonometric degree of
the generating functions $\chi_0^{(r)}$, $\chi_1^{(r)}$, $X_2^{(r)}$
and $Y_2^{(r)}$ is equal to $rK\,$, the choice to set $K=2$ and the
rule~(c3) explain why we stopped the algorithm after having ended the
normalization step with $r=9\,$.

The behavior of the norms of the generating functions is reported in
fig.~\ref{fig:norm_gen_fun}.  Let us make a few comments.  The
theoretical estimates assure the convergence of the normal form
provided the norms decrease geometrically with the order.  The figure
shows that this is indeed the behavior in our case.  The behavior of
the $\Dscr_2^{(r)}$ functions appears to be quite irregular, but we
should recall that these functions do not play the role of {\it
  normalization\/} function, since they represent the {\it
  diagonalization\/} of the quadratic part related to the secular
variables.  We emphasize that the presence of a dangerous resonance
would be reflected in a sudden increase of the coefficients; thus, the
plot may be considered as a practical confirmation that the
frequencies are well chosen.

Concerning the computational effort, performing the construction of
the normal form up to order $r=9$ has taken about $10$~hours of CPU
time on an Intel Core i7, using about $6$~Gb of RAM.

\subsection{Explicit Calculation of the Orbits on the Elliptic Torus}
\label{sbs:expl_calc}

We now perform a check on the approximation of the elliptic torus.  To
this end, we calculate the orbit on the torus using the analytic
expression and we compare it with a numerical integration of
Hamilton's equations.  In this subsection we explain how the
calculation of the orbit via normal form is performed.

According to the theory of Lie series, the canonical transformation
$(\vet{p},\vet{q},\vet{x},\vet{y})=
\Kscr^{(r)}\big(\vet{p}^{(r)},\vet{q}^{(r)},\vet{x}^{(r)},\vet{y}^{(r)}\big)$
inducing the normalization up to the step $r$ is given by
\begin{equation}
\vcenter{\openup1\jot\halign{
 \hbox {\hfil $\displaystyle {#}$}
&\hbox {\hfil $\displaystyle {#}$\hfil}
&\hbox {$\displaystyle {#}$\hfil}\cr
\Kscr^{(r)}\big(\vet{p}^{(r)},\vet{q}^{(r)},\vet{x}^{(r)},\vet{y}^{(r)}\big)
&=
&\exp\lie{\Dscr_2^{(r)}}\circ\,\exp\lie{Y_2^{(r)}}\circ\,\exp\lie{X_2^{(r)}}
\circ\,
\cr
& &\exp\lie{\chi_1^{(r)}}\circ\,\exp\lie{\chi_0^{(r)}}\circ\,\ldots\,\circ\,
\exp\lie{\Dscr_2^{(1)}}\circ\,\exp\lie{Y_2^{(1)}}\circ\,
\cr
& &\exp\lie{X_2^{(1)}}
\circ\,\exp\lie{\chi_1^{(1)}}\circ\,\exp\lie{\chi_0^{(1)}}
\>\big(\vet{p}^{(r)},\vet{q}^{(r)},\vet{x}^{(r)},\vet{y}^{(r)}\big)\,,
\cr
}}
\label{eq:trasf_Kolmogorov-like}
\end{equation}
where
$\big(\vet{p}^{(r)},\vet{q}^{(r)},\vet{x}^{(r)},\vet{y}^{(r)}\big)$
are meant to be the new coordinates. Thus, the canonical
transformation $(\vet{p},\vet{q},\vet{x},\vet{y})=
\Kscr^{(\infty)}(\vet{P},\vet{Q},\vet{X},\vet{Y})$ brings the initial
Hamiltonian $H^{(0)}$ in the normal form
$H^{(\infty)}=H^{(0)}\circ\Kscr^{(\infty)}$, which is written
in~(\ref{Ham^(infty)}), with
$\Kscr^{(\infty)}=\lim_{r\to\infty}\Kscr^{(r)}$. Let us introduce a
new symbol to denote the composition of all the canonical changes of
coordinates defined in sects.~\ref{sec:2D_plan_Ham}
and~\ref{sec:form_alg}, i.e.,
\begin{equation}
\Cscr^{(r)}=\Escr\circ\Tscr_{\vet{\Lambda}^*}\circ\Dscr\circ\Kscr^{(r)}\,,
\label{def-Cscr}
\end{equation}
where $(\tilde{\vet{r}},\vet{r})=
\Escr(\vet{\Lambda},\vet{\lambda},\vet{\csi},\vet{\eta})$ is the
canonical transformation giving the heliocentric positions $\vet{r}$
and their conjugated momenta $\tilde{\vet{r}}$ as a function of the
Poincar\'e's variables. If $\big(\tilde{\vet{r}}(0),\vet{r}(0)\big)$
is an initial condition on an invariant elliptic torus, in principle
we might use the following calculation scheme to integrate the
equation of motion:
\begin{equation}
\vcenter{\openup1\jot\halign{
 \hbox to 12 ex{\hfil $\displaystyle {#}$\hfil}
&\hbox to 12 ex{\hfil $\displaystyle {#}$\hfil}
&\hbox to 54 ex{\hfil $\displaystyle {#}$\hfil}\cr
\big(\tilde{\vet{r}}(0),\vet{r}(0)\big)
&\build{\longrightarrow}_{}^{{{\displaystyle
\left(\Cscr^{(\infty)}\right)^{-1}}
\atop \phantom{0}}}
&\left({{\displaystyle \vet{P}(0)=\vet{0}}
\,,\, {\displaystyle \vet{Q}(0)}\,,\,{\displaystyle \vet{X}(0)=\vet{0}}
\,,\, {\displaystyle \vet{Y}(0)}=\vet{0}}\right)
\cr
& &\bigg\downarrow \build{\Phi_{\vet{\omega}^{(\infty)}\cdot\vet{P}}^{t}}_{}^{}
\cr
\big(\tilde{\vet{r}}(t),\vet{r}(t)\big)
&\build{\longleftarrow}_{}^{{{\displaystyle 
\Cscr^{(\infty)}} \atop \phantom{0}}}
&\left({{\displaystyle  \vet{P}(t)=\vet{0}}
\,,\, {\displaystyle \vet{Q}(t)=\vet{Q}(0)+\vet{\omega}^{(\infty)} t}
\,,\,{\displaystyle \vet{X}(t)=\vet{0}}
\,,\, {\displaystyle \vet{Y}(t)}=\vet{0}}\right)
\cr
}}
\ \ \,,
\label{semi-analytical_scheme}
\end{equation}
where $\Phi_{\vet{\omega}^{(\infty)}\cdot\vet{P}}^{t}$ induces the
quasi-periodic flow related to the frequencies vector
$\vet{\omega}^{(\infty)}$.  Of course, the previous scheme requires an
unlimited computing power; from a practical point of view, we can just
approximate it by replacing $\Cscr^{(\infty)}$ with $\Cscr^{(R)}$,
being $R$ as large as possible.  Thus, the integration via normal form
actually reduces to a transformation of the coordinates of the initial
point to the coordinates of the normal form, the calculation of the
flow at time $t$ in the latter coordinates, which is a trivial matter
since the flow is exactly quasi-periodic with known frequencies,
followed by a transformation back to the original coordinates.

The rough but natural check that one can perform consists precisely in
comparing the orbit calculated via the semi-analytic scheme just
described with the result of a direct numerical integration.  As it
has been shown in~\cite{Loc-Gio-2000}, \cite{Loc-Gio-2005},
\cite{Loc-Gio-2007} and~\cite{Gab-Jor-Loc-2005}, this kind of
comparisons provide a very stressing test for the accuracy of the
whole algorithm constructing the normal form.

\begin{figure}
\centerline{\includegraphics[width=175mm]{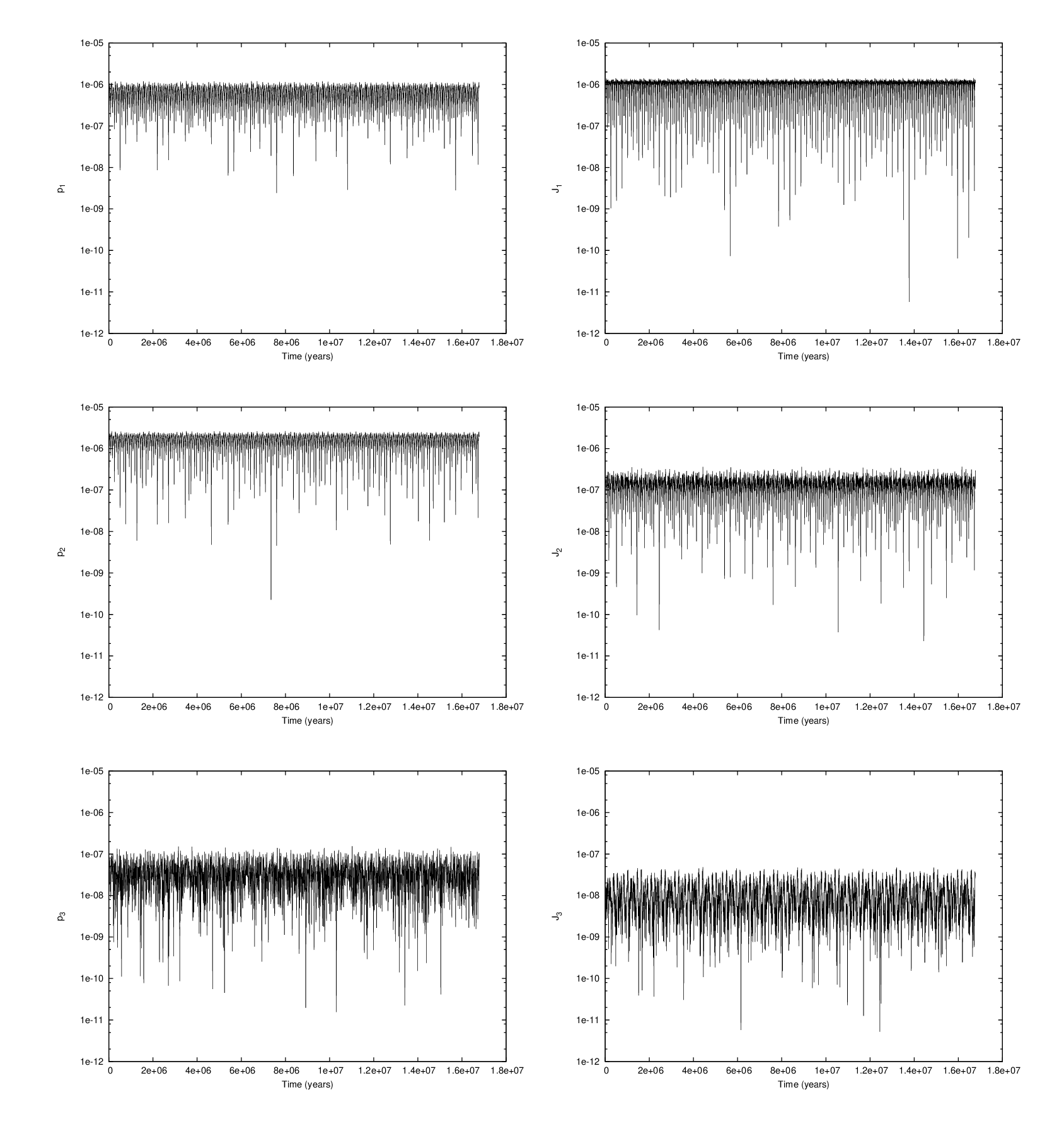}}
\caption{Time evolution of the normalized fast actions
  $p_1,\,p_2,\,p_3$ and of the slow actions $J_1,\,J_2,\,J_3$ along a
  numerically integrated orbit.  These quantities should be zero on
  the elliptic torus, and one sees here that they take indeed very
  small values along the orbit, thus showing that the approximation of
  the elliptic torus is good enough.}
\label{fig:errors}
\end{figure}

The initial condition $\big(\tilde{\vet{r}}(0),\vet{r}(0)\big)$ of the
orbit that we want to study is not directly provided by the
observations.  Thus, here we perform a test in a slightly different
way with respect to previous papers.  Precisely, we check that the
normalized variables $\vet{P}\,$, $\vet{X}$ and $\vet{Y}$ remain
constantly equal to zero along the motion on an elliptic torus.  Let
$p_1,\, p_2,\, p_3$ and $J_1,\,J_1,\, J_3$ denote the fast actions and
the secular actions of the Hamiltonian normalized up to order $R$ (we
omit the upper index $R$ for simplicity).  Setting $p_1(0) = p_2(0) =
p_3(0) = 0$ and $J_1(0) = J_2(0) = J_3(0) = 0$ we have a point on the
elliptic torus, and by applying the transformations we get the initial
condition $\big(\tilde{\vet{r}}(0),\vet{r}(0)\big)\,$.  More
precisely, we focus on the case $R=9$ and we consider the initial
conditions
\begin{equation}
\left(\Cscr^{(9)}\right)\left(\vet{0},\vet{0},\vet{0},\vet{0}\right)\,;
\label{cond-in-toro-ell-dopo-9-passi}
\end{equation}
according to the calculation method described in this section.  This
should be an accurate approximation of a point on an elliptic torus.
Therefore, we preliminarily integrated the motion of the planar SJSU
system over a time interval of $2^{24}$ years, by using the symplectic
method $\Sscr\Bscr\Ascr\Bscr_{3}$ (see~\cite{Las-Rob-2001}) with a
time-step of~$0.04$ years.  For a set of points on the orbit we apply
the canonical transformation $\big(\Cscr^{(9)}\big)^{-1}$, and so we
calculate the evolved actions $p_1(t)\,$, $p_2(t)\,$, $p_3(t)$ and
$J_1(t)\,$, $J_2(t)\,$, $J_3(t)\,$, which should be zero for all $t$
apart from a small oscillation due to the approximation of the
invariant torus. The results are illustrated in fig.~\ref{fig:errors},
where it is seen that the value of the actions remains very small,
reaching about $10^{-6}$ in the worst cases.

\subsection{Validation of the Results by Using Frequency
  Analysis}\label{sbs:val_freq_anal}

Another more refined test is based on frequency map analysis~(see,
e.g.,~\cite{Laskar-99} and~\cite{Laskar-2005} for an introduction).
The key point is that, according to the ideal calculation
scheme~(\ref{semi-analytical_scheme}), the time evolution of every
pair of canonical coordinates (e.g. the maps $t\mapsto\tilde
r_l(t)+\imunit r_l(t)$ with $l=1$, $2$, $3\,$) may be expressed as
\begin{equation}
\sum_{j=0}^{\infty}c_j\exp\big(\imunit\zeta_jt\big)\,,
\quad
{\rm where},\ \forall\ j\ge 0\,,\,\ c_j\in\complessi\ {\rm and}
\ \,\exists\ \,\vet{k}_j\in\interi^{n_1}
\ {\rm such\ that\ }\zeta_j=\vet{k}_j\cdot\vet{\omega}^{(\infty)}\,.
\label{spettro-quasi-periodico}
\end{equation}
Thus, the frequency spectrum is a very peculiar one, since it contains
only combinations of the fast frequencies, while no secular
frequencies should appear.  By the way, the same peculiarity of the
spectrum is characteristic of any signal produced by calculating the
evolution of a pair of mechanical quantities along an orbit on the
elliptic torus.  This is what we want to check in the present section.
From a strictly mathematical point of view, let us recall that the
previous formula for the Fourier spectrum can be deduced by the
scheme~(\ref{semi-analytical_scheme}), because the conjugacy function
$\vet{Q}\mapsto\Cscr^{(\infty)}\left(\vet{0},\vet{Q},\vet{0},\vet{0}\right)$
is analytic, as follows from the convergence of the constructive
algorithm that we shall prove in a forthcoming paper.

The orbit with initial point as in the previous
sect.~\ref{sbs:expl_calc} has been sampled with a time interval of 1
year.\footnote{Here we should add a remark concerning the precision of
  the calculation.  In order to have a signal clean enough to be
  analyzed a particular care about the precision is mandatory.  After
  some trials tuning the parameters of the numerical integration, we
  found that the 80 bits floating point numbers provided by the
  current AMD and INTEL CPUs fits our needs.  Technically this is
  obtained by using the {\tt long double} types of the GNU {\bf C}
  compiler under a Linux operating system.}  The fast frequencies
$\vet{\omega}^{(\infty)}$ have been accurately determined by looking
at the main components of the Fourier spectrum of the signals
$\Lambda_l(t)\exp\big({\rm i}\lambda_l(t)\big)\,$, with
$l=1,\,2,\,3\,$.  As test functions we used the time evolution of the
secular Poincar\'e's variables, namely $\csi_l(t)+\imunit\eta_l(t)$
with $l=1$, $2$, $3\,$.  The corresponding signals have been submitted
to the frequency analysis method using the so-called Hanning filter.
The reason for the choice of the coordinates as test functions is that
if the initial point given by~(\ref{semi-analytical_scheme}) is not on
the elliptic invariant torus then the presence of secular frequencies
is expected to be particularly evident in these signals.

\begin{table*}
\caption[]{Decomposition of the Fourier spectrum of the signal
  $\csi_3(t)+\imunit\eta_3(t)\,$, which is related to the Uranus
  secular motion. The numerical values in the table have been obtained
  by applying the frequency analysis method. See the text for more
  details.}
\label{tab:scomp_spettr_Urano}
\begin{center}
  \begin{tabular}{|c|c|c|c|c|}
\hline
$j$ &$\zeta_j$ &$\vet{k}_j$
&$|\zeta_j-\vet{k}_j\cdot\vet{\omega}^{(\infty)}|$ &$|c_j|$
\\
\hline
0 &$-7.48019221455542005\times\,10^{-2}$
&$(0,0,-1)$ &$0.0\times\,10^{+00}$ &$2.9770\times\,10^{-4}$
\\
1  &$\phantom{-}3.80127210702886631\times\,10^{-1}$
&$(1,0,-2)$  &$5.6\times\,10^{-17}$ &$5.5428\times\,10^{-5}$
\\
2  &$\phantom{-}6.37064849761184715\times\,10^{-2}$
&$(0,1,-2)$  &$5.6\times\,10^{-17}$ &$1.8199\times\,10^{-5}$
\\
3  &$\phantom{-}2.02214892097791255\times\,10^{-1}$
&$(0,2,-3)$  &$0.0\times\,10^{+00}$ &$1.7410\times\,10^{-5}$
\\
4  &$\phantom{-}3.40723299219463982\times\,10^{-1}$
&$(0,3,-4)$  &$0.0\times\,10^{+00}$ &$6.1013\times\,10^{-6}$
\\
5  &$\phantom{-}8.35056343551327518\times\,10^{-1}$
&$(2,0,-3)$  &$5.6\times\,10^{-17}$ &$3.7452\times\,10^{-6}$
\\
6  &$\phantom{-}4.79231706341136876\times\,10^{-1}$
&$(0,4,-5)$  &$1.7\times\,10^{-16}$ &$2.4485\times\,10^{-6}$
\\
7 &$-2.13310329267227178\times\,10^{-1}$
&$(0,-1,0)$  &$2.5\times\,10^{-16}$ &$1.4521\times\,10^{-6}$
\\
8  &$-3.51818736388899878\times\,10^{-1}$
&$(0,-2,1)$  &$2.2\times\,10^{-16}$ &$1.0175\times\,10^{-6}$
\\
9 &$\phantom{-}6.17740113462809326\times\,10^{-1}$
&$(0,5,-6)$  &$1.1\times\,10^{-16}$ &$1.0447\times\,10^{-6}$
\\
10  &$\phantom{-}1.28998547639976824\times\,10^{+0}$
&$(3,0,-4)$  &$2.2\times\,10^{-16}$ &$7.8098\times\,10^{-7}$
\\
11 &$-4.90327143510572494\times\,10^{-1}$
&$(0,-3,2)$  &$1.1\times\,10^{-16}$ &$7.1175\times\,10^{-7}$
\\
12  &$\phantom{-}7.56248520584482442\times\,10^{-1}$
&$(0,6,-7)$  &$2.2\times\,10^{-16}$ &$4.6141\times\,10^{-7}$
\\
13 &$-9.84660187842435808\times\,10^{-1}$
&$(-2,0,1)$  &$1.7\times\,10^{-16}$ &$4.1885\times\,10^{-7}$
\\
14 &$-6.28835550632244611\times\,10^{-1}$
&$(0,-4,3)$  &$5.0\times\,10^{-16}$ &$3.8157\times\,10^{-7}$
\\
15 &$-5.29731054993994532\times\,10^{-1}$
&$(-1,0,0)$  &$5.6\times\,10^{-16}$ &$2.9840\times\,10^{-7}$
\\
16 &$-1.11363990387936461\times\,10^{-5}$
&$(0,0,0)$  &$1.1\times\,10^{-05}$ &$2.2654\times\,10^{-7}$
\\
17 &$\phantom{-}8.94756927706155003\times\,10^{-1}$
&$(0,7,-8)$  &$1.1\times\,10^{-16}$ &$2.0832\times\,10^{-7}$
\\
18  &$-7.67343957753918837\times\,10^{-1}$
&$(0,-5,4)$  &$1.0\times\,10^{-15}$ &$1.9160\times\,10^{-7}$
\\
19  &$\phantom{-}1.74491460924820907\times\,10^{+0}$
&$(4,0,-5)$  &$2.8\times\,10^{-16}$ &$1.7777\times\,10^{-7}$
\\
20 &$-1.43958932069087675\times\,10^{+0}$
&$(-3,0,2)$  &$1.1\times\,10^{-16}$ &$1.3450\times\,10^{-7}$
\\
21  &$\phantom{-}2.43025734926846093\times\,10^{-2}$
&$(-1,4,-4)$ &$1.1\times\,10^{-14}$ &$1.1086\times\,10^{-7}$
\\
22  &$-9.05852364875589622\times\,10^{-1}$
&$(0,-6,5)$  &$9.4\times\,10^{-16}$ &$9.3984\times\,10^{-8}$
\\
23 &$\phantom{-}1.03326533482782756\times\,10^{+0}$
&$(0,8,-9)$  &$0.0\times\,10^{+00}$ &$9.5496\times\,10^{-8}$
\\
24  &$-1.96924221667578817\times\,10^{-5}$
&$(0,0,0)$  &$2.0\times\,10^{-05}$ &$5.2900\times\,10^{-8}$
\\
\hline
\end{tabular}
\end{center}
\end{table*}

In table~\ref{tab:scomp_spettr_Urano} we report our numerical results
about the first $25$ terms of the
decomposition~(\ref{spettro-quasi-periodico}) for the Uranus secular
signal, i.e., $\csi_3(t)+\imunit\eta_3(t)\,$.  The vectors
$\vet{k}_j\in\interi^{n_1}$ listed in the third column are determined
so that the absolute difference
$|\zeta_j-\vet{k}_j\cdot\vet{\omega}^{(\infty)}|$ is minimized, with
the limit $|\vet{k}_j|\le 20\,$.  We stress that some limits on the
absolute value of $\vet{k}_j$ must be imposed in order to make
consistent its calculation, and our choice is motivated by the fact
that the decay of Fourier coefficients in the analytic conjugacy
function
$\Cscr^{(\infty)}\left(\vet{0},\vet{Q},\vet{0},\vet{0}\right)$ leads
us to expect that the main contributions to the spectrum are related
to low order harmonics.

If the initial conditions~(\ref{cond-in-toro-ell-dopo-9-passi}) were
exactly on an elliptic torus, each value
$|\zeta_j-\vet{k}_j\cdot\vet{\omega}^{(\infty)}|$ reported in the
fourth column of table~\ref{tab:scomp_spettr_Urano} should be equal to
zero.  We see that all of them, except for the cases corresponding to
$j=16,\,24\,$, are actually small enough to be considered as generated
by round-off errors.  On the other hand, we can definitely say that
$\zeta_{16}\simeq -1.1\times 10^{-5}$ and $\zeta_{24}\simeq -1.9\times
10^{-5}$ are ``secular frequencies'', because their values are
$\Oscr(\mu)\,$.  Indeed, let us recall that $\mu\simeq 10^{-3}$, but
the mass ratio for Uranus, i.e., $m_3/m_0\simeq 4.4\times 10^{-5}$, is
even smaller.

Let us say that a total absence of secular frequencies can occur only
in a very ideal situation, namely: (i)~all the calculations described
in sects.~\ref{sec:2D_plan_Ham} and~\ref{sec:form_alg} should be
carried out without performing any truncations on the expansions,
(ii)~the initial conditions~(\ref{cond-in-toro-ell-dopo-9-passi})
should be replaced with $\Cscr^{(\infty)}
\left(\vet{0},\vet{0},\vet{0},\vet{0}\right)\,$, (iii)~no numerical
errors should be there.  Concerning the point~(iii), the numerical
error does clearly appear in our signal, but it turns out to be very
small.  The points~(i) and~(ii) are more relevant, because our
analytical calculation gives a point close to, but not lying on the
elliptic torus.  Thus the presence of secular frequencies should be
expected, and the size of the corresponding coefficients may be
considered as a measure of the distance of the initial point from the
invariant torus.  Nevertheless, it is very remarkable that the
amplitude of the first found secular frequency, i.e., $|c_{16}|\,$, is
three orders of magnitude smaller than the main component of the
spectrum.  We consider this fact as a clear indication that our
algorithm is able to identify the elliptic torus with a remarkable
precision.

Other components corresponding to secular frequencies are expected to
be even smaller than those found with $j=16,\,24\,$.  In fact, let us
recall that the frequency analysis method detects the summands
$c_j\exp\big(\imunit\zeta_jt\big)$ appearing
in~(\ref{spettro-quasi-periodico}) in a nearly decreasing order with
respect to the amplitude $|c_j|$ (for instance, just two exchanges are
needed in order to rewrite table~\ref{tab:scomp_spettr_Urano} in the
correct decreasing order); moreover, we calculated that the
discrepancy $\big|\csi_3(t)+\imunit\eta_3(t)-
\sum_{j=0}^{24}c_j\exp\big(\imunit\zeta_jt\big)\big|$ is smaller than
about $\simeq 3.7\times 10^{-7}$ for all the time values~$t$ for which
we sampled the signal.  Let us emphasize that such an upper bound on
the maximal discrepancy is just a little larger than the amplitude
$|c_{16}|\,$.

A similar decomposition has been calculated for both the signals
$\csi_1(t)+\imunit\eta_1(t)$ and $\csi_2(t)+\imunit\eta_2(t)$ (which
are related to the secular motions of Jupiter and Saturn,
respectively).  The behavior is very similar to that of
table~\ref{tab:scomp_spettr_Urano}, so we omit the corresponding
tables.

In order to check the effectiveness of our algorithm, let us further
investigate the results as a function of the normalization step in a few
different ways.

\begin{figure}
\centerline{\includegraphics[width=120mm]{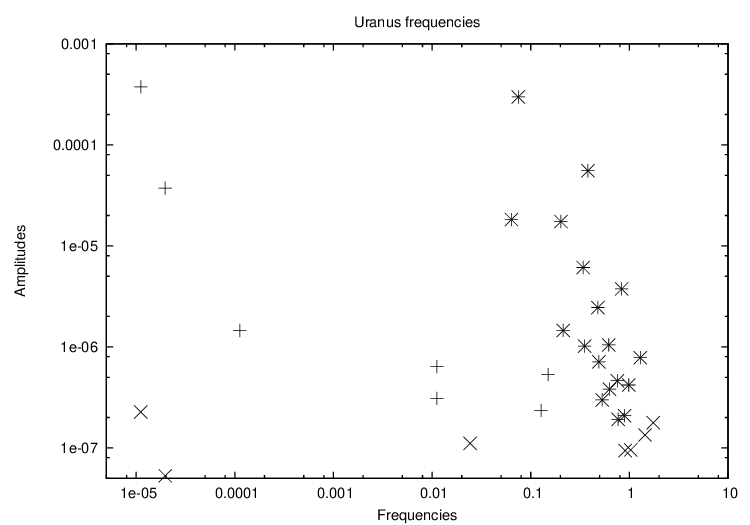}}
\caption{Frequency analysis of the secular signal related to the
  secular Uranus motion: $\csi_3(t)+\imunit\eta_3(t) =
  \sum_{j=0}^{\infty}c_j\exp\big(\imunit\zeta_jt\big)\,$. Plot of the
  amplitudes $|c_j|$ as a function of the frequencies $\zeta_j$ in
  Log-Log scale. The symbol $\times$ [$+$, resp.]  refer to the signal
  related to the motion starting from the initial
  conditions~(\ref{cond-in-toro-ell-dopo-9-passi})
  [(\ref{cond-in-toro-ell-dopo-0-passi}), resp.], i.e., the
  approximation of a point on an elliptic torus after having
  performed~$9$ [$0$, resp.]  steps of the algorithm constructing the
  corresponding normal form.  In both cases, the results for just the
  first~$25$ components have been reported in the figure above.}
\label{fig:ainf_Urano}
\end{figure}

First, in fig.~\ref{fig:ainf_Urano} we plot the amplitudes $|c_j|$
as a function of the frequencies $\zeta_j$ in Log-Log scale.  We
compare the amplitudes corresponding to two different orbits, namely:
(a)~the orbit with initial point calculated as
in~(\ref{cond-in-toro-ell-dopo-9-passi}), and (b)~the orbit with initial point 
\begin{equation}
\Cscr^{(0)}\left(\vet{0},\vet{0},\vet{0},\vet{0}\right)\,.
\label{cond-in-toro-ell-dopo-0-passi}
\end{equation}

The symbol $\times$ and $+$ correspond to the orbits (a)~and (b),
respectively.  So the data marked with $\times$ refer to the
decomposition of table~\ref{tab:scomp_spettr_Urano}.  Let us remark
that $\Cscr^{(0)}\left(\vet{0},\vet{0},\vet{0},\vet{0}\right)=
\Escr\circ\Tscr_{\vet{\Lambda}^*}\circ\Dscr
\left(\vet{0},\vet{0},\vet{0},\vet{0}\right)$ is a sort of trivial
approximation of a point on the elliptic torus as it is provided by
simply avoiding to apply the part of our algorithm constructing the
normal form, as it is described in sect.~\ref{sec:form_alg}.  Let us
recall that the amplitudes related to secular frequencies are expected
to have a size $\Oscr(\mu)\,$.  Thus, they are easily separated from
the fast ones, having a size less than $10^{-3}$.  By looking at the
right side of fig.~\ref{fig:ainf_Urano}, one immediately sees that the
part of the spectrum related to the fast frequencies is nearly
indistinguishable when the initial
conditions~(\ref{cond-in-toro-ell-dopo-9-passi})
or~(\ref{cond-in-toro-ell-dopo-0-passi}) are considered, because the
symbols $\times$ and $+$ superpose each other almost exactly for most
of the main components. On the other hand, the secular part of the
spectrum, that is in the left side of fig.~\ref{fig:ainf_Urano},
strongly differs. In fact, when the initial
conditions~(\ref{cond-in-toro-ell-dopo-0-passi}) (that trivially
approximate a point on the elliptic torus) are considered, three
secular frequencies are detected; while just two of them are found in
the case of the more accurate initial
data~(\ref{cond-in-toro-ell-dopo-9-passi}). Moreover, by comparing the
amplitudes, one can see that the better approximation of the torus
with the orbit~(a) makes the corresponding secular amplitudes
($\times$) to be approximately three orders of magnitude smaller than the
amplitudes ($+$) for orbit~(b).  The results are very similar for the
three planets, thus we reported only those for Uranus.

\begin{figure}
\centerline{\includegraphics[width=120mm]{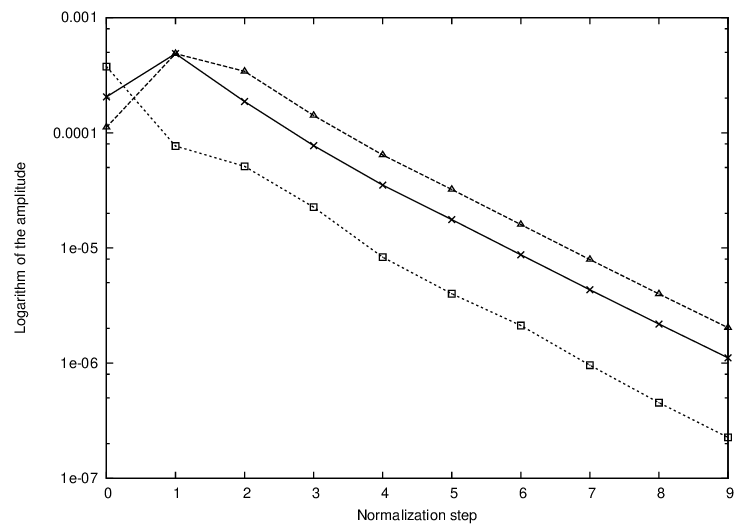}}
\caption{Frequency analysis of the signal related to the secular
  motion of Jupiter, Saturn and Uranus, i.e., $\csi(t)+\imunit\eta(t)
  = \sum_{j=0}^{\infty}c_j\exp\big(\imunit\zeta_jt\big)\,$, where
  $\csi$ and $\eta$ are the secular Poincar\'e variables of a planet
  (recall definition~\eqref{var-Poincare-piano}). Plot of the
  logarithm of the amplitude $|c_j|$ related to the main secular
  component as a function of the normalization step. The symbols
  $\times\,$, $\triangle\,$ and $\square\,$ refer to signals
  corresponding to Jupiter, Saturn and Uranus, respectively.}
\label{fig:decadono}
\end{figure}

For completeness we report in fig.~\ref{fig:decadono} the amplitude of
the main secular component of the signals for the three planets as a
function of the normalization step.  The symbols $\times\,$,
$\triangle\,$ and $\square\,$ refer to the data for Jupiter, Saturn
and Uranus, respectively.  A sharp geometric decay to zero is
observed, which makes evident the effectiveness of our procedure
constructing the normal form for an elliptic torus.

Finally, it is interesting to look at the change of the secular
frequencies during the normalization procedure.  The values of the
secular frequencies $\vet{\Omega}^{(r)}$ in~\eqref{Ham^(r)} of the
three planets, calculated at every step of the normalization procedure
are collected in the table~\ref{tab:sec}.  At step $0$ the frequencies
so obtained correspond to those given by Lagrange-Laplace theory for
our planar model including Sun, Jupiter, Saturn and Uranus.  After the
first few steps a steady convergence to a final value is highlighted.

\begin{table*}
\caption[]{Values of the frequencies $\vet{\Omega}^{(r)}$ of the
  transversal oscillation with respect to the elliptic torus.}
\begin{center}
\label{tab:sec}
\begin{tabular}{|c|c|c|c|}
\hline
Step ($r$) & ${\Omega}^{(r)}_0$ & ${\Omega}^{(r)}_1$ & ${\Omega}^{(r)}_2$ \\
\hline
$0$ & $-1.0779281\times\,10^{-4}$ & $-1.7700018\times\,10^{-5}$ & $-1.1005216\times\,10^{-5}$ \\
\hline
$1$ & $-1.0837401\times\,10^{-4}$ & $-1.7718250\times\,10^{-5}$ & $-1.0998189\times\,10^{-5}$ \\
\hline
$2$ & $-1.1413488\times\,10^{-4}$ & $-1.9540195\times\,10^{-5}$ & $-1.1137012\times\,10^{-5}$ \\
\hline
$3$ & $-1.1211306\times\,10^{-4}$ & $-1.9684669\times\,10^{-5}$ & $-1.1133050\times\,10^{-5}$ \\
\hline
$4$ & $-1.1207431\times\,10^{-4}$ & $-1.9686809\times\,10^{-5}$ & $-1.1133969\times\,10^{-5}$ \\
\hline
$5$ & $-1.1207181\times\,10^{-4}$ & $-1.9691417\times\,10^{-5}$ & $-1.1135412\times\,10^{-5}$ \\
\hline
$6$ & $-1.1206555\times\,10^{-4}$ & $-1.9691851\times\,10^{-5}$ & $-1.1135625\times\,10^{-5}$ \\
\hline
$7$ & $-1.1206559\times\,10^{-4}$ & $-1.9691857\times\,10^{-5}$ & $-1.1135808\times\,10^{-5}$ \\
\hline
$8$ & $-1.1206555\times\,10^{-4}$ & $-1.9691864\times\,10^{-5}$ & $-1.1135815\times\,10^{-5}$ \\
\hline
$9$ & $-1.1206555\times\,10^{-4}$ & $-1.9691864\times\,10^{-5}$ & $-1.1135815\times\,10^{-5}$ \\
\hline
\end{tabular}
\end{center}
\end{table*}

\subsection*{Acknowledgments}
The authors have been supported by the research program ``Dynamical
Systems and applications'', PRIN 2007B3RBEY, financed by MIUR.
M.S. has been partially supported also by the research program ``Studi
di Esplorazione del Sistema Solare'', financed by ASI.



\begin{thebibliography}{50}
\bibitem{Arnold-1963.1}{V.I.~Arnold: \emph{ Proof of a theorem of A. N.
    Kolmogorov on the invariance of quasi-periodic motions under
    small perturbations of the Hamiltonian}, Usp. Mat. Nauk, {\bf 18},
  13 (1963); Russ. Math. Surv., {\bf 18}, 9 (1963).}
\bibitem{Arnold-1963.2}{V.I.~Arnold: \emph{ Small denominators and
    problems of stability of motion in classical and celestial
    mechanics}, Usp. Math. Nauk {\bf 18} N.6, 91 (1963);
  Russ. Math. Surv. {\bf 18} N.6, 85 (1963).}
\bibitem{Ben-Gal-Gio-Str-1984}{G.~Benettin, L.~Galgani, A.~Giorgilli and
  J.M.~Strelcyn: \emph{A Proof of Kolmogorov's Theorem on Invariant
    Tori Using Canonical Transformations Defined by the Lie method},
  Nuovo Cimento, {\bf 79}, 201--223 (1984).}
\bibitem{Bia-Chi-Val-2003}{L.~Biasco, L.~Chierchia and E.~Valdinoci:
  \emph{Elliptic two-dimensional invariant tori for the planetary
    three-body problem}, Arch. Rational Mech. Anal., {\bf 170} ,
  91--135 (2003).}
\bibitem{Bia-Chi-Val-2006}{L.~Biasco, L.~Chierchia and E.~Valdinoci:
  \emph{N-dimensional elliptic invariant tori for the planar
    (N+1)-body problem}, SIAM Journal
  on Mathematical Analysis, {\bf 37} , n. 5, 1560--1588 (2006).}
\bibitem{Birkhoff-27}{G.D.~Birkhoff: \emph{Dynamical systems}, New
  York (1927).}
\bibitem{Cel-Gio-Loc-2000}{A.~Celletti, A.~Giorgilli and U.~Locatelli:
  \emph{ Improved Estimates on the Existence of Invariant Tori for
    Hamiltonian Systems}, Nonlinearity, {\bf 13}, 397--412 (2000).}
\bibitem{Cas-Jor-2000}{E.~Castell\`a and A.~Jorba: \emph{On the
    vertical families of two-dimensional tori near the triangular
    points of the Bicircular problem}, Cel. Mech. \& Dyn. Astr., {\bf
    76}, 35--54 (2000).}
\bibitem{Deprit-1983}{A.~Deprit: \emph{Elimination of the nodes in
    problems of $n$ bodies}, Cel. Mech. \& Dyn. Astr., {\bf 30},
  181--195 (1983).}
\bibitem{Gab-Jor-2001}{F.~Gabern and A.~Jorba: \emph{A restricted
    four-body model for the dynamics near the Lagrangian points of the
    Sun--Jupiter system}, DCDS-B, {\bf 1}, 143--182 (2001).}
\bibitem{Gab-Jor-Loc-2005}{F.~Gabern, A.~Jorba and U.~Locatelli:
  \emph{ On the construction of the Kolmogorov normal form for the
    Trojan asteroids}, Nonlinearity, {\bf 18}, n.4, 1705--1734
  (2005).}
\bibitem{Giorgilli-1995}{A.~Giorgilli: \emph{Quantitative methods in
  classical perturbation theory}, proceedings of the Nato ASI school
  ``From Newton to chaos: modern techniques for understanding and
  coping with chaos in N-body dynamical systems'', A.E.\ Roy e
  B.D.\ Steves eds., Plenum Press, New York (1995).}
\bibitem{Gio-Loc-1997.1}{A.~Giorgilli and U.~Locatelli:
  \emph{Kolmogorov theorem and classical perturbation theory}, J. of
  App. Math. and Phys. (ZAMP), {\bf 48}, 220--261 (1997).}
\bibitem{Gio-Loc-1997.2}{A.~Giorgilli and U.~Locatelli: \emph{On
    classical series expansion for quasi-periodic motions}, MPEJ, {\bf
    3}, 5, 1--25 (1997).}
\bibitem{Gio-Loc-San-2009}{A.~Giorgilli, U.~Locatelli and
  M.~Sansottera: \emph{Kolmogorov and Nekhoroshev theory for the
    problem of three bodies}, Cel. Mech. \& Dyn. Astr., {\bf 104},
  159--173 (2009).}
\bibitem{Gio-Loc-San-2010}{A.~Giorgilli, U.~Locatelli and
  M.~Sansottera: \emph{Su un'estensione della teoria di Lagrange per i
    moti secolari}, Istituto Lombardo (Rend. Scienze), {\bf 143},
  223--240 (2009).}
\bibitem{Gio-Loc-San-2011}{A.~Giorgilli, U.~Locatelli and
  M.~Sansottera: \emph{A constructive algorithm of the normal form for
    lower dimensional elliptic tori.}, in preparation.}
\bibitem{Jef-Mos-1966}{W.H.~Jefferys and J.~Moser:
  \emph{Quasi-periodic solutions for the three-body problem},
  Astronom.~J., {\bf 71}, 568--578 (1966).}
\bibitem{Jor-Vil-1997.1}{A.~Jorba and J.~Villanueva: \emph{On the
    persistence of lower dimensional invariant tori under
    quasiperiodic perturbations}, J. of Nonlin. Sci., {\bf 7}, 427--473
  (1997).}
\bibitem{Jor-Vil-1997.2}{A.~Jorba and J.~Villanueva: \emph{On the
    Normal Behaviour of Partially Elliptic Lower Dimensional Tori of
    Hamiltonian Systems}, Nonlinearity, {\bf 10}, 783--822 (1997).}
\bibitem{Jor-Vil-1998}{A.~Jorba and J.~Villanueva: \emph{Numerical
    Computation of Normal Forms Around Some Periodic Orbits of the
    Restricted Three Body Problem}, Physica D, {\bf 114}, 197--229
  (1998).}
\bibitem{Kolmogorov-1954}{A.N.~Kolmogorov: \emph{ Preservation of
    conditionally periodic movements with small change in the Hamilton
    function}, Dokl. Akad. Nauk SSSR, {\bf 98}, 527
  (1954). Engl. transl. in: Los Alamos Scientific Laboratory
  translation LA-TR-71-67; reprinted in: Lecture Notes in Physics {\bf
    93}.}
\bibitem{Lagrange-1776}{J.L.~Lagrange: \emph{Sur l'alt\'eration des
  moyens mouvements des plan\`etes}, Mem. Acad. Sci. Berlin {\bf 199}
  (1776); \emph{Oeuvres compl\`etes}, {\bf VI}, 255, Paris,
  Gauthier--Villars (1869).}
\bibitem{Lagrange-1781}{J.L.~Lagrange: \emph{Th\'eorie des variations
  s\'eculaires des \'el\'ements des plan\`etes. Premi\`ere partie
  contenant les principes et les formules g\'en\'erales pour
  d\'e\-ter\-mi\-ner ces variations}, Nouveaux m\'emoires de
  l\:Acad\'emie des Sciences et Belles-Lettres de Berlin (1781);
  \emph{Oeuvres compl\`etes}, {\bf V}, 125--207, Paris, Gauthier--Villars
  (1870).}
\bibitem{Lagrange-1782}{J.L.~Lagrange: \emph{Th\'eorie des variations
  s\'eculaires des \'el\'ements des plan\`etes. Seconde partie
  contenant la d\'etermination de ces variations pour chacune des
  plan\`etes pricipales}, Nouveaux m\'emoires de l\:Acad\'emie des
  Sciences et Belles-Lettres de Berlin (1782); \emph{Oeuvres
  compl\`etes}, {\bf V}, 211--489, Paris, Gauthier--Villars (1870).}
\bibitem{Laplace-1772}{P.S.~Laplace: \emph{M\'emoire sur les solutions
  particuli\`eres des \'equations diff\'erentielles et sur les
  in\'egalit\'es s\'eculaires des plan\`etes} (1772); \emph{Oeuvres
  compl\`etes}, {\bf IX}, 325, Paris, Gauthier-Villars (1895).}
\bibitem{Laplace-1784}{P.S.~Laplace: \emph{M\'emoire sur les
  in\'egalit\'es s\'eculaires des plan\`etes et des satellites},
  Mem. Acad. royale des Sci. de Paris (1784); \emph{Oeuvres
    compl\`etes}, {\bf XI}, 49, Paris, Gauthier-Villars (1895).}
\bibitem{Laplace-1785}{P.S.~Laplace: \emph{Th\'eorie de Jupiter et de
  Saturne}, Mem. Acad. royale des Sci. de Paris (1785); \emph{Oeuvres
  compl\`etes}, {\bf XI}, 164, Paris, Gauthier-Villars (1895).}
\bibitem{Laskar-1989b}{J.~Laskar: \emph{Syst\`emes de variables et
  \'el\'ements}, in Benest, D. and Froeschl\'e, C. (eds.): \emph{Les
  M\'ethodes modernes de la M\'ecanique C\'eleste}, 63--87, Editions
  Fronti\`eres (1989).}
\bibitem{Las-Rob-1995}{J.~Laskar and P.~Robutel: \emph{Stability of
  the Planetary Three-Body Problem --- I. Expansion of the Planetary
  Hamiltonian}, Celestial Mechanics and Dynamical Astronomy, {\bf 62},
  193--217 (1995).}
\bibitem{Laskar-99}{J.~Laskar: \emph{Introduction to frequency map
  analysis}, in C. Sim\`o (managing ed.), Proceedings of the NATO ASI
  school: ``Hamiltonian Systems with Three or More Degrees of
  Freedom'', S'Agaro (Spain), June 19--30, 1995, Kluwer, 134--150
  (1999).}
\bibitem{Laskar-2005}{J.~Laskar: \emph{Frequency Map analysis and quasi
  periodic decompositions}, in Benest et al. (managing eds):
  ``Hamiltonian systems and Fourier analysis'', Taylor and Francis
  (2005).}
\bibitem{Las-Rob-2001}{J.~Laskar and P.~Robutel: \emph{High order
  symplectic integrators for perturbed {H}amiltonian systems},
  Celestial Mechanics and Dynamical Astronomy, {\bf 80}, 39--62
  (2001).}
\bibitem{Lieberman-1971}{B.~Lieberman: \emph{Existence of
    quasi-periodic solutions to the three-body problem}, Cel. Mech. \&
  Dyn. Astr., {\bf 3}, 408--426 (1971).}
\bibitem{Loc-Gio-2000}{U.~Locatelli and A.~Giorgilli: \emph{Invariant
    tori in the secular motions of the three-body planetary systems},
  Cel. Mech. \& Dyn. Astr., {\bf 78}, 47--74 (2000).}
\bibitem{Loc-Gio-2005}{U.~Locatelli and A.~Giorgilli:
  \emph{Construction of the Kolmogorov's normal form for a planetary
    system}, Regular and Chaotic Dynamics, {\bf 10}, n.2, 153--171
  (2005).}
\bibitem{Loc-Gio-2007}{U.~Locatelli and A.~Giorgilli: \emph{Invariant
    tori in the Sun--Jupiter--Saturn system}, DCDS-B, {\bf 7},
  377--398 (2007).}
\bibitem{Mal-Rob-Las-02}{F.~Malige, P.~Robutel and J.~Laskar:
  \emph{Partial reduction in the N-body planetary problem using the
    angular momentum integral}, Celestial Mechanics and Dynamical
  Astronomy, {\bf 84}, 283--316 (2002).}
\bibitem{Moser-1962}{J.~Moser: \emph{On invariant curves of
    area-preserving mappings of an annulus},
  Nachr. Akad. Wiss. G\"ott,. II Math. Phys. Kl 1962, 1--20 (1962).}
\bibitem{Pinzari-tesi}{G.~Pinzari: \emph{On the Kolmogorov set for
    many-body problems}, Ph.D. thesis, Universit\`a di Roma Tre
  (2009); publicly available at the web page: {\tt
    http://ricerca.mat.uniroma3.it/dottorato/Tesi/pinzari.pdf}\thinspace.}
\bibitem{Poincare-1892}{H.~Poincar\'e: \emph{Les m\'ethodes nouvelles
    de la M\'ecanique C\'eleste}, Gauthier-Villars, Paris (1892),
  reprinted by Blanchard (1987).}
\bibitem{Poincare-1905}{H.~Poincar\'e: \emph{Le\c cons de M\'ecanique
    C\'eleste}, tomes I--II, Gauthier-Villars, Paris (1905).}
\bibitem{Poschel-1989}{J.~P{\"o}schel: \emph{On elliptic lower
  dimensional tori in Hamiltonian systems}, Math. Z., {\bf 202},
  559--608 (1989).}
\bibitem{Poschel-1996}{J.~P{\"o}schel: \emph{A KAM-theorem for some
  nonlinear PDEs}, Ann. Scuola Norm. Pisa Cl. Sci., {\bf 23}, 119--148
  (1996).}
\bibitem{Robutel-1995}{P.~Robutel: \emph{Stability of the Planetary
    Three-Body Problem --- II. KAM Theory and Existence of
    Quasiperiodic Motions}, Celestial Mechanics and Dynamical
  Astronomy, {\bf 62}, 219--261 (1995).}
\bibitem{San-Loc-Gio-2010}{M.~Sansottera, U.~Locatelli and
  A.~Giorgilli: \emph{On the stability of the secular evolution of the
    planar Sun--Jupiter--Saturn--Uranus system}, Math. Comput. Simul.,
  {\tt doi:10.1016/j.matcom.2010.11.018} (2011).}
\bibitem{Standish-1998}{E.M.~Standish: \emph{JPL Planetary and Lunar
    Ephemerides, DE405/LE405}, Jet Propulsion Laboratory ---
  Interoffice memorandum, {\bf IOM 312.F -- 98 -- 048} (1998).}
\end{thebibliography}
\end{document}